\outer\def\give#1. {\medbreak
             \noindent{\bf#1. }}                     
\outer\def\section #1\par{\bigbreak\centerline{\S
     {\bf#1}}\nobreak\smallskip\noindent}
\def\({\left(}
\def\){\right)}
\def\contract{{\,{\vrule height5pt width0.4pt depth 0pt}
{\vrule height0.4pt width6pt depth 0pt}\,}}

\def\sqr#1#2{{\vcenter{\hrule height.#2pt              
     \hbox{\vrule width.#2pt height#1pt\kern#1pt
     \vrule width.#2pt}
     \hrule height.#2pt}}}
\def\square{\mathchoice\sqr{5.5}4\sqr{5.0}4\sqr{4.8}3\sqr{4.8}3}
\def\qed{\hskip4pt plus1fill\ $\square$\par\medbreak}



\def\cA{{\cal A}}

\def\cC{{\cal C}}
\def\cD{{\cal D}}

\def\cG{{\cal G}}

\def\cL{{\cal L}}

\def\cQ{{\cal Q}}
\def\cR{{\cal R}}

\def\cV{{\cal V}}
\def\cW{{\cal W}}

 

\def\C{{\bf C}}
\def\cx#1{{\C}^{#1}}     
\def\cp1{{{\bf P}^1}}


\def\bar{\overline}              
 
 


\magnification=\magstep1

\centerline{Polynomial Diffeomorphisms of $\cx 2$:}
\centerline{V. Critical Points and Lyapunov Exponents}
\bigskip
\centerline{Eric Bedford and John Smillie}
\bigskip
\section 0. Introduction

This paper deals with the dynamics of polynomial diffeomorphisms
$f:\C^2\to\C^2$. To exclude trivial cases we make the standing assumption that
the dynamical degree $d=d(f)$  is greater than
one (see Section 1 for a definition).  It is often quite useful in dynamics to
focus attention on invariant objects. A natural invariant set to consider is
$K=K_f$, the set of points with bounded orbits. Pluripotential theory allows us
to associate to this set the harmonic measure, $\mu=\mu_f$ of $K_f$. 
For polynomial diffeomorphisms this measure is finite and invariant, and we
normalize it to have total mass one.  In previous papers we have shown that this
measure has considerable dynamical significance. We have shown that $\mu$ is
ergodic [BS3] and that the support of $\mu$ is the closure of the set of periodic
saddle orbits [BLS1]. Further, $\mu$ is the unique measure of maximal entropy
[BLS1], and $\mu$ describes the distribution of periodic points [BLS2].

To any measure we can associate Lyapunov exponents.  The rate of expansion and
contraction of tangent vectors at a point $p$  by $f$ is measured by a pair of
Lyapunov exponents, $\lambda^+(p)$ and $\lambda^-(p)$. In the presence of an
ergodic invariant measure such as $\mu$ these exponents are constant almost
everywhere, and we denote them by $\lambda^+(\mu)$ and $\lambda^-(\mu)$.  By
[BS3] the (complex) Lyapunov exponents of $\mu$
satisfy $\lambda^-(\mu)<0<\lambda^+(\mu)$.  This condition is known as
(nonuniform) hyperbolicity of the measure $\mu$.  Nonuniform hyperbolicity
implies that at $\mu$ almost every point $p$ there is a spitting of the tangent
space into  complex one dimensional subspaces $E^u(p)$ and $E^s(p)$ so that  for
$v\in E^u(p)$ we have
$||Df^n(v)||\sim \exp(n\lambda^+)||v||$ and for $v\in E^s(p)$ we have
$||Df^n(v)||\sim \exp(n\lambda^-)||v||$.
In this paper we will prove an integral formula for the Lyapunov exponents.  In
many ways our formula is analogous to the Brolin-Manning formula for Lyapunov
exponents with respect to harmonic measure for polynomial maps of ${\bf C}$,
which we now describe.

Let $g$ be a polynomial map of $\C$. We
let $K_g$ denote the set of points with bounded orbits. We denote by $\mu=\mu_g$
the harmonic measure of $K_g$.  There is a single Lyapunov exponent
$\lambda(\mu)$ which gives the average rate of expansion  along the orbit $\mu$
almost everywhere.  The Green function of $K$ is given by the following formula:
$$G(z)=\lim_{n\to\infty}d^{-n}\log^+|g^n(z)|,$$
which relates it to the superexponential rate at which an orbit escapes to
infinity.  The following Brolin-Manning formula relates
the Lyapunov exponents to the critical $c_j$ points of the map:
$$\lambda(\mu)=\log d+\sum G(c_j).\eqno(0.1)$$
 The above formula was obtained in the case without critical points
by Manning [M]; the present formulation appears in Przytycki [Pr].

 Formula (0.1) takes an especially simple form in the quadratic case,
$g(z)=z^2+c$. If we write $G_c$ for the Green function of the Julia
set of $g(z)=z^2+c$, we have
$$\lambda(\mu)=\log 2+G_c(0).$$
Douady and Hubbard [DH] observed that $G_c(0)={1\over2}h(c)$, where $h(c)$ is
the Green function of the Mandelbrot set. Thus the Lyapunov exponent is
connected to the potential theory of the parameter space.  The idea of
understanding the parameter space by means of potential theory sparked our
interest in Lyapunov exponents.  Some properties of the function
$\Lambda(f)=\lambda^+(\mu)$ as a function on the  parameter space of polynomial
diffeomorphisms were studied in [BS3].

In this paper we do two things.  First we define a notion of critical point and
critical point measure for  polynomial diffeomorphisms of $\C^2$ and explore the
dynamical significance of these objects. Second, we use this measure to prove
an integral formula which is the analog of (0.1).

One ingredient in our integral formula is a Green function. The
function $G$ has two analogs in $\cx2$:
$$G^\pm(x,y)=\lim_{n\to\infty}d^{-n}\log^+|f^{\pm n}(x,y)|.$$
We write $K^\pm\subset\cx2$ for the set of points in $\cx2$ bounded in
forward/backward time and we define $U^\pm$ to be the complement of $K^\pm$. The 
functions $G^\pm$ are zero on $K^\pm$ and pluriharmonic on $U^\pm$ and serve as
Green functions. The function $G^+$ which describes the forward rate of escape
is the analog of $G$.

The function $f$ is a diffeomorphism and hence has no critical points in the
usual sense of the word. For maps in one variable critical points of $g$ with
unbounded orbits are associated to critical points of the Green function $G$.
This suggests that we look for critical points of $G^+$. Since $\nabla G^+$ is
non-zero at every point of $U^+$, $G^+$ has no critical points in the usual
sense. In many situations the best analog of the set $\C$ for polynomial maps is
not all of $\C^2$ but is rather the set $K^-$ of points with bounded backward
orbits. This suggests that we should look for critical points of $G^+$
restricted to $K^-$. We could make sense of this concept if $K^-$ were a
manifold. Now $K^-$ is a rather wild set and in particular it is not a manifold.
On the other hand for $\mu$ almost every point in $J$ the set $W^u(p)$ is a
manifold and this manifold is contained in $J^-\subset K^-$. We define our
set of unstable critical points, $\cC^u$, to be the set of critical points of the
restrictions $G^+|W^u(p)$.

These critical points have an interesting dynamical interpretation which does not
make explicit reference to the function $G^+$. In the region $U^+$ points
escape to infinity in forward time. In fact they escape at a super-exponential
rate. In $U^+$ there is a plane field $\tau^+$ such that for $v\in\tau^+$,
$Df^n(v)$ decreases super-exponentially as $n\to\pm\infty$ (Lemma 1.2).  We will
show that a critical point in $W^u(p)$ as defined above is a point at which the
tangent space of $W^u(p)$ concides with $\tau^+$. In [HO] a holomorphic foliation
$\cG^\pm$ of $U^\pm$ was constructed. We will show in Proposition B.1 that the
tangent space to a leaf of this foliaton is given by $\tau^\pm$, and the global
leaves of $\cG^\pm$ are super-stable/-unstable manifolds. Thus critical points
are points at which unstable manifolds $W^u(p)$ and super-stable manifolds
intersect tangentially and we can think of such a critical point as a type of
heteroclinic tangency.

The next ingredient in our integral formula is a critical measure $\mu_c^-$
supported on the set of critical points. In order to define this measure, we use
the unstable current $\mu^-={1\over2\pi}dd^cG^+$.  The unstable current $\mu^-$,
in some sense, serves as the current of integration on the unstable lamination
$\cW^u$.  Locally, $\mu^-$ may be thought of as follows.  There is a
(pairwise disjoint) family of unstable disks $D_t$ and a transversal measure on
the space of parameters $t$; and $\mu^-$ is given locally as the current of
integration over
$D_t$, integrated with respect to the transverse measure.  We construct the
critical measure $\mu^-_c$ by replacing the current of integration over $D_t$ by
the sum of the point masses at the critical points of $G^+|D_t$. 

Another way to approach the critical measure is to fix (arbitrarily) a vector
$\alpha$ and a covector $\beta$.  The variety
$$Z_k(\alpha,\beta)=\{x:\beta\cdot Df^k_x(\alpha)=0\}$$
is the set of points $x$ where $Df^k_x$ maps the vector $\alpha$ to the kernel
of $\beta$.  If $M\subset U^+$ is a Riemann surface such that $G^+|_M$ is not
locally constant, then by Lemma 5.2, $Z_k(\alpha,\beta)\cap M$ converges to
the set of critical points of $G^+|_M$ as $k\to\infty$.  The slice of the current
$\mu^-$ by the variety $f^jZ_k(\alpha,\beta)$ is given by the wedge product
$\mu^-\wedge [f^jZ_k(\alpha,\beta)]$.  We show (Theorem 5.9) that the average
(over
$\alpha$) of these slices gives the critical measure:
$$\mu^-_c=\lim_{j\to\infty\atop k-j\to\infty}\int\sigma(\alpha)\,
\mu^-\wedge[f^jZ_k(\alpha,\beta)].$$ 
The set $\cC^s$ of critical points in $J^+$ and the measure $\mu^+_c$ can be
defined in an analogous way.

The main result in this paper (see Theorem 6.1 and Corollary 6.6) is a formula
for the Lyapunov exponents of harmonic measure:
\proclaim Theorem.
$$\eqalign{\lambda^+(\mu)&=\log d+\int_{\{1\le
G^+<d\}}G^+\mu^-_c\cr
\lambda^-(\mu)&=-\log d-\int_{\{1\le
G^-<d\}}G^-\mu^+_c.\cr}\eqno(0.2)$$

\noindent (We omit the conventional ``$d$'' from in front of the measure in the
integral in order to reduce confusion with the exterior derivative operator and
with the degree.)  The condition $\{1\le G^+<d\}$ in the formula for $\lambda^+$
has the effect of choosing a fundamental domain for the action of
$f$ on the set $\cC^u$. This insures that each orbit of critical points
contributes only once. Other choices of fundamental domains work equally well.
This is a consequence of the fact that the integrand $G^+\mu^-$ is invariant
under $f$. The invariance of the integrand occurs because
$G^+$ multiplies by a factor of $d$, and $\mu^-$ multiplies by a factor of
$d^{-1}$ under $f$.  A geometrically appealing way of finding a fundamental
domain arises naturally when $f$ generates a real horseshoe mapping, and in this
case the formula may be given in a particularly simple form; see Appendix A. 

We will briefly describe some of the connections between Lyapunov exponents,
dimension of harmonic measure and connectivity for polynomial maps and
polynomial diffeomorphisms. 

The Hausdorff dimension of a measure is, by definition, the
infimum of the Hausdorff dimensions of Borel sets with full measure.  Let $J$
be the Julia set of a polynomial map $g$ of the complex plane.  The Lyapunov
exponent of $g$  with respect to harmonic measure $\mu_J$ is related to the
Hausdorff dimension of $\mu_J$ by the formula $\lambda(\mu)\hbox{\rm
HD}(\mu)=\log d$.  Formula (0.1) has the consequence that the Hausdorff
dimension of the harmonic measure of the Julia set is at most one and is equal
to one if and only if all critical points have bounded orbits.  Thus the Julia
set is connected if and only if the harmonic measure has Hausdorff dimension
one. (It is a general result of Makarov that the Hausdorff dimension of the
harmonic measure of any connected set is one.)

For polynomial diffeomorphisms there is also a connection between exponents and
certain planar sets.  Let $f$ be a polynomial diffeomorphisms of $\cx2$ and let
$\mu$ denote the corresponding harmonic measure.  For $\mu$ almost every point
$p$ stable and unstable manifolds $W^{s/u}(p)$ exist and are complex manifolds
conformally equivalent to $\cx{}$.  Given such a $p$ we can consider the sets
$W^u(p)\cap K^+$ and $W^s(p)\cap K^-$ (which can be viewed as subsets of the
complex plance.)   The ``slice measures'' $\mu^\pm|_{W^{s/u}(p)}$ play the role
of harmonic measures for the sets $K^\pm\subset W^{u/s}$.  In [BLS1] the slice
measures were shown to satisfy the Ledrappier-Young [LY] formula
$$\lambda^+(\mu)={\log d\over \hbox{HD}(\mu^-|_{W^u_{loc}(p)}) }, \qquad
\lambda^-(\mu)={\log d\over \hbox{HD}(\mu^+|_{W^s_{loc}(p)}) }.
$$
We explore the relation between exponents, critical points and
connectivity for polynomial diffeomorphisms in [BS6].

Critical points play an important role in the dynamical study of polynomial
maps. We have defined two sets of ``critical points" namely $\cC^s$ and $\cC^u$,
but there are other possible definitions that could be made. The critical points
in $\cC^s$ and $\cC^u$ are points at which  there is a vector $v$ with the
property that $Df^n(v)$ decreases in both forward and backward time.  If we take
this condition to be a characteristic of critical points we can also ask about
the set $\cC$ of ``critical points'' for which both the forward and backward
orbits are unbounded.  These are points in $U^+\cap U^-$ at which the
super-stable and super-unstable foliations are tangent. (Such points of tangency
were first considered by Hubbard.)  For a given polynomial diffeomorphism $f$
either of the sets
$\cC^s$ or $\cC^u$ may be empty, but we show in
Proposition B.3 that the set $\cC$ is never empty. In this paper we do not
discuss the remaining case of ``critical points''  with bounded
forward and backward orbits, which is more difficult (see [BC] and [BY]).

A question that arises for diffeomorphisms of $\cx2$ is the relation between
stable and unstable critical points.  The integral formula allows us to deduce
(Corollary 6.9) that if $f$ is dissipative, then $\cC^s\ne\emptyset$, which by
[BS6] is seen to have topological consequences.

Section 1 contains introductory material and an analysis of the growth of
tangent vectors in $U^+$. In \S2 we develop the laminarity
of $\mu^+$ using elementary methods.  We show that for any algebraic variety
$X$, the convergence of $cd^{-n}[f^nX]$ to $\mu^-$ induces the geometry of
the laminar structure.  
In \S3 we introduce results from Pesin theory
concerning the stable/unstable manifolds with respect to the hyperbolic
measure $\mu$.  We show (Lemma 3.3) that the laminarity of the convergence of
$f^nX$ also respects the laminar structure of the Pesin manifolds.
In \S4 we define alternative notions of average rates of growth and we relate
these to the Lyapunov exponent.
In \S5 we introduce unstable critical measure $\mu^-_c$ and establish some of
its properties. In \S6 we prove the integral formula.

\section 1.  Super-stable directions in $U^+$

The polynomial diffeomorphisms of $\cx2$ were classified up to conjugacy by
Friedland and Milnor [FM], who introduced a degree $d$, which we call dynamical
degree, and which is defined by the formula $d=\lim_{n\to\infty}({\rm
deg}(f^n))^{1/n}$.  They showed that the mappings with interesting dynamics
correspond to the case  $d\ge2$; and so, as in the one-dimensional case, we will
consider only polynomial diffeomorphisms with $d\ge2$.

The inverse of a polynomial diffeomorphism is again a polynomial diffeomorphism,
and for this reason polynomial diffeomorphisms form a group and are often called
polynomial automorphisms. Friedland and Milnor showed that a dynamically
nontrivial mapping is conjugate in the group of polynomial diffeomorphisms to a
mapping of the form $f=f_1\circ\cdots\circ f_m$, where $f_j(x,y)=(y,p_j(y)-a_jx)$
and $p_j(y)=y^{d_j}+O(y^{d_j-2})$, with $d_j\ge2$.  The algebraic degree of $f$
is then $d=d_1\cdots d_m$, and the iterates for $k\ge1$ are given by
$$f^k(x,y)=(y^{d^k/d_1}+\cdots,y^{d^k}+\cdots).\eqno(1.2)$$
Thus for maps in this standard form the algebraic degree coincides with the
dynamic degree. In particular it is an integer greater than or equal to 2. 

Certain dynamical properties of $f$ can be deduced simply from a consideration
of the sizes of the coordinates.  We recall the following standard notations
and results. 
$$V^+=\{|y|\ge|x|,|y|\ge R\},\quad V^-=\{|y|\le|x|,|x|\ge R\},\quad
V=\{|x|,|y|\le R\}.\eqno(1.3)$$
There is an $R$ sufficiently large that $fV^+\subset V^+$, $fV\subset V\cup
V^+$, and that for any point $(x,y)\in V^-$ the orbit $f^n(x,y)$ can remain in
$V^-$ only for finitely many $n>0$.  Further, if $K^+$ is the set of points with
bounded forward orbits, then
$$U^+=\cx2-K^+=\bigcup_{n=0}^\infty f^{-n}V^+,\quad{\rm and\ }\quad
U^-=\cx2-K^-=\bigcup_{n=0}^\infty f^{n}V^-.\eqno(1.4)$$

The rate of escape to infinity in forward/backward time:
$$G^\pm=\lim_{n\to\infty}{1\over d^n}\log^+|f^{\pm n}(x,y)|\eqno(1.5)$$
links potential theory and dynamics.  The function $G^\pm$ is continuous and
pluri-subharmonic on $\cx2$.  We let $\pi_1$ and $\pi_2$ denote projections onto
the first and second coordinates.  So for a point $(x,y)\in V^+$, we have
$\pi_2f^n\sim(\pi_1f^n)^{d_1}$ as $n\to+\infty$.  Thus $|f^n|\sim|\pi_2f^n|$ as
$n\to+\infty$, and it follows that
$$G^+=\lim_{n\to\infty}{1\over d^n}\log^+|\pi_2f^{ n}(x,y)|.$$

In the remainder of this
section we will analyze the growth rates of tangent vectors at points in $U^+$. 
For $(x,y)\in V^+$,
$$|\pi_2f_1(x,y)-y^{d_1}|=O(y^{d_1-2})+O(x),$$
where the $O$ terms are uniform on $V^+$.  Thus the following estimate is
uniform in $n\ge0$, and $(x,y)\in V^+$:
$$\eqalign{ |\pi_2f^{n+1}-(\pi_2f^n)^d| &\le O((\pi_2f_2\circ\cdots\circ
f_m\circ f^n)^{d_1-2}) +O(\pi_1f_2\circ\cdots\circ f_m\circ f^n)\cr
&\le O((\pi_2f^n)^{(d_1-2)d_2\cdots d_m})+O((\pi_2f^n)^{d_2\cdots d_m})\cr
&\le O((\pi_2 f^n)^{d-2}).\cr}$$
For $(x,y)\in V^+$ and any $N$, we have
$$ G^+(x,y) ={1\over d^N}\log|\pi_2f^N| +
  \sum_{n=N}^\infty \left({1\over d^{n+1}}\log|\pi_2f^{n+1}|  -{1\over
d^{n}}\log|\pi_2f^{n}| \right)$$
By the estimates above, the $n$-th term in the series is dominated by
$$\eqalign{
{1\over d^{n+1}}\log\left|{\pi_2f^{n+1}\over (\pi_2f^n)^d}\right|&\le
{1\over d^{n+1}}
\log\left(1+\left|{\pi_2f^{n+1}-(\pi_2f^n))^d\over(\pi_2f^n)^d}\right|\right)\cr
&\le{1\over d^{n+1}}\log\left(1+{O(|\pi_2f^n|^{d-2})\over
|\pi_2f^n|^d}\right)\cr
&\le{C\over d^{n+1}|\pi_2f^n|^2}.\cr}$$
Summing the tail of the series, we conclude that
$$\left|G^+(x,y)-{1\over d^N}\log|\pi_2f^N(x,y)|\right|\le{C'\over
d^N|\pi_2f^N(x,y)|^2}\eqno(1.6)$$
holds for all $N$ and $(x,y)\in V^+$.

\give Remark.  The rate of convergence of numerical approximations to $G^+$ and
$\partial G^+$ can be understood by (1.6) and the estimates that precede it. 
Suppose that $f(x,y)=(y,y^d+q(y)-ax)$, where $\deg q\le d-2$.  We set
$y_{-1}=x$, $y_0=y$, and $y_{n+1}=y_n^d+q(y_n)-ay_{n-1}$ for $n\ge0$.  We may
approximate $G^+$ either by $d^{-N}\log|y_N|$ or by the telescoping sum
$$\log|y_k|+\sum_{n=k}^Nd^{-n-1}\log\left\vert {y_{n+1}/ y_n^d}\right\vert
=\log|y_k|+\sum_{n=k}^Nd^{-n-1}\log|1+\rho_n|,$$
where $\rho_n=(q(y_n)-ay_{n-1})y_n^{-d}=O(y_n^{-2})$.

Similarly, we set $\partial y_{-1}=dx$, $\partial y_0=dy$, and $\partial
y_{n+1}=(d\cdot y_n^{d-1}+q'(y_n))\partial y_n -a\partial y_{n-1}$ for $n\ge0$. 
Thus we may approximate $\partial G^+$ by $(2d^Ny_N)^{-1}\partial y_N$ or by the
telescoping sum
$${(2d^ky_k)^{-1}\partial y_k } + \sum_{n=k}^N{d^{-n}}\left({(2d\cdot
y_{n+1})^{-1}\partial y_{n+1}}-{(2y_n)^{-1}\partial y_n\ }\right),$$
which, after cancellation, is
$${\partial y_k\over 2d^ky_k} + \sum_{n=k}^N{d^{-n}(1+\rho_n)^{-1}}
\left( {-\rho_n\partial y_n\over 2y_n} +{q'(y_n)\partial y_n \over 2d\cdot
y_n^d}-{a\partial y_{n-1}\over 2d\cdot y_n^d}\right)$$
so the $n$-th term in the summation is no larger than
$O(d^{-n}y_n^{-2})$. \bigskip

For a tangent vector $v$, we will use the notation $\partial G\cdot v$ for the
pairing with the 1-form $\partial G$, and $Df(v)$ for the action of the
differential $Df$.

\proclaim Lemma 1.1.  There exist $C$ and $R$ sufficiently large that
$$\left|\partial G^+\cdot v -{1\over d^n}{\partial(\pi_2f^n)\cdot
v\over|\pi_2f^n|}\right| \le {C|v|\over d^n|f^n|^2}$$
holds for all $n\ge0$, all $(x,y)\in V^+$, and all tangent vectors $v\in
T_{(x,y)}\cx2$.

\give Proof.  We have estimate (1.6) in a neighborhood of fixed radius about
any point of $V^+$.  Since these functions are harmonic, we may differentiate
this estimate and have the same estimate also for the gradients.
\qed

The following gives a dichotomy on the growth rate of $Df^n(v)$.  Either it
grows super-exponentially to $\infty$ as $n\to\infty$, or it decreases to 0
super-exponentially.

\proclaim Lemma 1.2.  If $(x,y)\in U^+$ and $v\in T_{(x,y)}\cx2$ is a vector 
with $\partial G^+\cdot v=0$, then there are constants $c$ and $N$ such that
$$|Df^n( v)|\le{c|v|\over|f^n|}.\eqno(1.7)$$
for $n\ge N$.  If $\partial G^+\cdot v\ne0$, then
$$|Df^n( v)|\sim d^n|f^n|\,|\partial G^+\cdot v|.\eqno(1.8)$$
Further, if $(x_0,y_0)\in U^+$ and $\epsilon>0$ are given, then there exist
small $\delta>0$ and large $N$ such that
$$|Df^n( v)|\ge\delta d^n|f^n|\,|\partial G^+\cdot v|\eqno(1.9)$$
holds in a $\delta$ neighborhood of $(x_0,y_0)$ for all $n\ge N$ and all tangent
vectors $v$ such that $|\partial G^+\cdot v|\ge\epsilon|v|$.

\give Proof. To make estimates, we may identify $Df^n$ with the
pair $(\partial(\pi_1 f^n),\partial(\pi_2 f^n))$.  By (1.1) and (1.6) it is
sufficient to estimate $\partial\pi_2 f^n$.  Thus (1.7) follows directly from
Lemma 1.1.  Again by Lemma 1.1, we estimate
$$|\partial(\pi_2 f^n)\cdot v|\ge d^n|\partial G^+\cdot v||\pi_2f^n|-{C|v|\over
|f^n|},$$
which yields (1.8) and (1.9).\qed

For $(x,y)\in U^\pm$ we let $\tau^\pm(x,y)$ denote the subspace of
$T_{(x,y)}\cx2$ annihilated by $\partial G^\pm$, i.e.\ such that $\partial
G^\pm\cdot v=0$ for all $v\in\tau^\pm$.  We will refer to $\tau^\pm$ as the
forward/backward dynamical critical directions.  If
$v\notin\tau^\pm$, then $|Df^n\cdot v|$ grows as $n\to\pm\infty$ at the rate
given in Lemma 1.2.

\proclaim Corollary 1.3.  If $v\in\tau^\pm$, then $\lim_{n\to\pm\infty}{1\over
|n|}\log|Df^n( v)|=-\infty$.

\give Proof.  By (1.7), 
$${1\over
|n|}\log|Df^n( v)|\le{1\over |n|}(\log|cv|-\log|f^n|)$$
which tends to $-\infty$ as $n\to\pm\infty$, since $|f^n|\sim
e^{d^{|n|}G^\pm}$. \qed

Next we consider another way to measure minimal growth of $Df^n$.  Let us fix
$n\in{\bf Z}$ and
$(x,y)\in\cx2$ and consider the mapping
$$T_{(x,y)}\cx2\ni t\mapsto{|Df^n\cdot t|\over|t|}.\eqno(1.10)$$
We let $\tau_n(x,y)$ denote the subspace of $T_{(x,y)}\cx2$ on which this
mapping is minimized.

\proclaim Proposition 1.4.  $\lim_{n\to+\infty}\tau_n=\tau^+$ on $U^+$ and
$\lim_{n\to-\infty}\tau_n=\tau^-$ on $U^-$.

\give Proof.   Let us suppose that there are vectors $v_{n_j}\in \tau_{n_j}$
which stay at positive angle from $\tau^+$.  Then there is an $\epsilon>0$ such
that $|\partial G^+\cdot v_{n_j}|\ge\epsilon|v_{n_j}|$.  By (1.9) of Lemma 1.2,
$|Df^{n_j}|$ grows as $j\to\infty$.  On the other hand, if $v^+\in\tau^+$, then
$Df^{n_j}v^+$ decreases to 0 as $j\to\infty$.  Thus for some large $j$,
$v_{n_j}$ does not minimize (1.10), which is a contradiction. \qed

\section 2. Laminar Properties of the Stable/Unstable Currents

In this section and the next we will discuss the laminar properties of the
currents $\mu^\pm$.   Laminarity is a ``natural'' structure for $\mu^\pm$ and has
been the key for understanding the deeper properties of $\mu^\pm$ and $\mu$.  It
will also be central to the definition of the critical measure.  In this section
we describe an explicit approach to laminarity which will be useful in section
5. In \S3 we describe an alternate approach to laminarity via the Pesin theory.
Although we will work only with $\mu^-$, it is evident that the analogous
properties hold for $\mu^+$. 

Let us summarize some notation and terminology about currents.  More details
are given in [BLS1].  We let $\cD_k$ denote the set of compactly supported
$k$-forms (test forms).  The dual space $\cD_k'$ is the set of $k$-dimensional
currents.  A sequence $\{T_n\}$ converges in the sense of currents if
$\lim_{n\to\infty} T_n(\varphi)=T(\varphi)$ for every test form
$\varphi\in\cD_k$.  If $X$ is a $k$-dimensional submanifold with locally finite
$k$-dimensional area, then there is the current of integration $[X]\in\cD_k'$,
whose action on a test form is given by 
$$[X](\varphi):=\int_X\varphi.$$
If $S$ is a discrete (0-dimensional) set, then the current of integration
$$[S]=\sum_{a\in S}\delta_a$$ is the sum of point masses at $S$.  It will be
useful for us to define the mass norm of a current as
$${\bf M}[T]=\sup_{|\varphi|\le1}|T(\varphi)|,$$
where $|\varphi|:=\sup_{x}|\varphi(x)|$ is the Euclidean supremum norm of a
test form $\varphi$.  The mass norm of $T$ is finite if and only if $T$ may be
represented as a linear combination of $k$-vectors with coefficients which are
finite measures.

If $\varphi$ is a smooth $k$-form, we define $T\contract\varphi$ by
$(T\contract\varphi)(\psi)=T(\varphi\wedge\psi)$.  If $S$ is a Borel set,
$T\contract S$ will denote the restriction of $T$ to $S$, i.e.
$T\contract\chi_S$, where $\chi_S$ is the characteristic function of $S$.  We
may do this whenever the mass norm of $T$ is locally finite.

While the stable/unstable currents $\mu^\pm:={1\over2\pi}dd^cG^\pm$ are
defined as positive, closed, (1,1)-currents,  they also
have special properties not enjoyed by general currents, and in fact it is
these properties that are the most useful for studying the dynamical
properties of $\mu^\pm$.  If $M$ is a 1-dimensional complex submanifold of
$\cx2$, then $\mu^\pm$  induce measures on $M$, given by
$$\mu^\pm|_M={1\over2\pi}(dd^c)|_M(G^\pm|_M)\eqno(2.0)$$
where $(dd^c)_M$ is the induced operator on $M$. 

If $\nu$ is a measure on the space $A$, and $\psi$ is an integrable function
on $A$, then we denote the integral of $\psi$ with respect to $\nu$ as
$\int\psi(a)\nu(a)$ or $\int\nu(a)\psi(a)$.  If $\{T_a:a\in A\}$ is a measurable
family of currents, we define the (direct) integral $\int\nu(a) T_a$ by its
action on a test form $\varphi$ by
$$\Big(\int_{a\in A}\nu(a)\, T_a\Big)(\varphi):=\int T_a(\varphi)\,\nu(a).$$
A current is laminar if it can be written as a direct integral of $T_a$ as
above, with the $T_a$ being currents of integration over pairwise disjoint
complex manifolds.

Our derivation of the laminar structure of $\mu^-$ will be based on the following
characterization of $\mu^-$.  Let $X\subset\cx2$ be an algebraic variety of pure
dimension 1.  By Proposition 4.2 of [BS1] there are positive integers $n_0$ and
$k$ such that $f^nX$ has degree $kd^{n-n_0}$ for $n\ge n_0$, and
$$\lim_{n\to\infty}{1\over kd^{n-n_0}}[f^nX]=\mu^-.\eqno(2.1)$$

All the constructions in this Section will depend on the variety $X$ and the
projection $\pi_\alpha(x,y)=\alpha_1x+\alpha_2y$ for some choice
of $\alpha\in\cx2$ with $|\alpha_1|^2+|\alpha_2|^2=1$.  We may rotate
corrdinates on $\cx2$ so that $\pi_\alpha(x,y)=y$, and we will not include
the choice of $X$ or $\alpha$ explicity in our notations.

Let $\cQ_s$ denote the set of squares in the plane with side $d^{-s}$ and with
vertices on points of the set $d^{-s}({\bf Z}+i{\bf Z})$.  Each square
$Q\in\cQ_s$ will be half-open, i.e.\ $Q=[a,b)\times[c,d)$, so that $\cQ_s$ is a
partition of $\cx{}$.   We choose $\kappa>0$ and let $\cQ_s'$ denote the set of
squares $Q'$ which have side of length $(1+2\kappa)d^{-s}$ and which are centered
about squares $Q$ of $\cQ_s$.  There is a number
$m(\kappa)$ such that each point of $\cx{}$ is contained in at most $m(\kappa)$
squares of $\cQ_s'$.  We let $Q_0$ denote a fixed square of $\cQ_0$.

Let $Q\subset Q_0$ be any square from $\cQ_s$, and let $Q'\in\cQ'_s$ denote the
square centered about it.  A connected component $\Gamma'$ of
$f^nX\cap\pi^{-1}Q'$ will be said to be {\it good} if the projection
$\pi|_{\Gamma'}:\Gamma'\to Q'$ is a homeomorphism.  We let 
$$\cG(Q,n)=\{\Gamma'\cap\pi^{-1}Q:\Gamma'{\rm\ good\ }\}.$$

Let us define
$$\mu^-_{\cQ_s,n}={1\over
kd^{n-n_0}}\sum_{Q\in\cQ_s}\sum_{\Gamma\in\cG(Q,n)}[\Gamma].$$
For fixed $Q_0\in\cQ_0$, there is a number $R>0$ such that
$$f^nX\cap\pi^{-1}Q'_0\subset\{|x|<R\}\eqno(2.2)$$
for all $n\ge n_0$.  

\proclaim Lemma 2.1.   There is a constant $C$,
independent of $Q\subset Q_0$, $s$ and $n$ such that  
$${\bf M}[(
k^{-1}d^{n_0-n}[f^nX]-\mu^-_{\cQ_s,n})\contract\pi^{-1}Q_0]\le
C{\rm Area}(Q)[1-k^{-1}d^{n_0-n}]m(\kappa).\eqno(2.3)$$

\give Proof.  The number of components of $f^nX\cap\pi^{-1}Q'$ is no more than
$kd^{n-n_0}$, the degree of $f^nX$.  If $\Gamma'$ is not good, then the number
of branch points, counted with multiplicity, in $\Gamma'$ is one less than the
mapping degree of $\pi|_{\Gamma'}$.  Thus the sum of the mapping degrees of
components that are not good is bounded above by $(kd^{n-n_0}-1)m(\kappa).$

Now we need to estimate $k^{-1}d^{n_0-n}$ times the area of the components that
are not good.  A property of analytic varieties (see Chirka [C]) is that there
is a constant $C$ depending on the $R$ of (2.2) and $\kappa$ such that the
area of every component $\Gamma'$ of $f^nX\cap\pi^{-1}Q'$ is bounded by
$${\rm Area}(\Gamma)\le C\mu\,{\rm Area}(Q)$$
where $\mu$ is the mapping degree of $\pi|_{\Gamma'}$.  Multiplying by
$d^{-2s}$, the area of $Q$, we estimate the mass in the left hand side of
(2.3) by the total of the mapping degrees coming from bad disks.  \qed

For each good $\Gamma'$ there is an analytic function $\varphi:Q'\to\cx{}$ such
that $\{(\varphi(y),y):y\in Q'\}=\Gamma'$.  Let $\cA(Q,n)$ denote the set of
all such analytic functions. 

Let us define 
$$S(Q,y,n)=\bigcup_{\varphi\in\cA(Q,n)}\{(\varphi(y),y)\}
=\bigcup_{\Gamma\in\cG(Q,n)}\Gamma\cap\pi^{-1}(y).
$$
The measures $\nu_Q(y,n):=k^{-1}d^{n_0-n}[S(Q,y,n)]$ are the slice measures of
$\mu^-_{Q,n}$ with respect to the projection $\pi$.   That is, $\nu_Q(y,n)$
is supported on $\pi^{-1}(y)$, and
$$\mu^-_{Q,n}\contract({i\over2}dy\wedge d\bar y) =\int_{y\in
Q}\cL^2(y)\,\nu_Q(y,n)$$
where $\cL^2$ denotes the Lebesgue area measure on $Q$.  Since the masses of the
currents $\mu^-_{Q,n}$ are uniformly bounded by (2.2), we may choose
a subsequence $\{n_j\}$ so that
$$\mu^-_Q:=\lim_{j\to\infty}\mu^-_{Q,n_j}\eqno(2.4)$$ 
exists. 
\proclaim Lemma 2.2. The limit
$\nu_Q(y):=\lim_{j\to\infty}\nu_Q({y,n_j})$ exists for every $y\in Q$.

\give Proof.  Each $\nu_Q(y,n)$ is a positive measure of mass at most one.  If
$\lim_{j\to\infty}\nu_Q(y,n_j)$ does not exist then there exist subsequences
$\{n^{(k)}_j\}$, $k=1,2$ of $\{n_j\}$ with distinct limiting measures
$\nu_Q^{(k)}=\lim_{j\to\infty}\nu_Q(y,n_j^{(k)})$, $k=1,2$.  We may assume that
there is a function $\varphi$ with 
$$\Big|\int\varphi(\nu^{1}_Q-\nu^{2}_Q)\Big|>\epsilon.$$

By the Cauchy estimate and (2.2), we have $|\psi'(y)|\le R/\kappa$ on $Q$ for
all $\psi\in\cG(Q,n)$.  Thus for any $n$,
$$\eqalign{
\Big|\int\varphi(\nu_Q(y_1,n)-\nu_Q(y_2,n))\Big|\le &
\sup\{|\varphi(a_1)-\varphi(a_2)|: \cr
&|a_1|,|a_2|\le R,
|a_1-a_2|\le{R\over\kappa}|y_1-y_2|\}.\cr}\eqno(2.5)$$
Now we choose $\delta$ such that $\sup\{|\varphi(a_1)-\varphi(a_2)|: |a_1|,
|a_2|\le R,|a_1-a_2|<R\delta/\kappa\}<\epsilon/2$.  Since the limit in (2.4)
exists, we have the limit $\nu_Q(\hat y)=\lim_{j\to\infty}\nu_Q(\hat y,n_j)$ for
almost every $\hat y\in Q$.  Thus we may choose $\hat y\in Q$ such that this
limit exists and such that $|\hat y-y|<\delta$.  It follows from (2.5) that
$$\Big|\int\varphi(\nu^{(k)}_Q-\nu_Q(\hat y,n_j))\Big|<\epsilon/2$$
for $j$ sufficiently large.  This is a contradiction, which proves the Lemma.
\qed

Let us set $S(Q):={\rm\ supp}\,\nu_Q(c_Q)$, and let $\cA(Q)$ denote the set of
analytic functions $\varphi:Q'\to\cx{}$ such that $\varphi(c_Q)\in S(Q)$, and
there is a sequence of functions $\varphi_{n_j}\in\cA(Q,n_j)$ converging to
$\varphi$. Since distinct good components must be disjoint, we have
$\varphi_1(y) \ne\varphi_2(y)$ for all $y\in Q'$, it follows that with the $R$ of
(2.2)
$$h(y)=\log(|2R|^{-1}|\varphi_1(y)-\varphi_2(y)|)$$
is a negative function on $Q'$.  By the Harnack inequality, there is a constant
independent of $s$, $n$, and $Q$ such that
$$h(y)\le const. \, h(c_Q)\quad{\rm for\ }y\in Q,$$
where $c_Q$ denotes the center of the square $Q$.  We conclude that there are
constants $R$ and $\kappa$ (independent of $s$, $n$, and $Q$) such that 
$$|\varphi_1(y)-\varphi_2(y)|\le
R|\varphi_1(c_Q)-\varphi_2(c_Q)|^\kappa\eqno(2.6)$$ 
for all $y\in Q$.

The following may be interpreted as a normal families argument for sets of
functions satisfying (2.6).

\proclaim Lemma 2.3. $\cA(Q)$ has the following properties:
\item{(1)}  For each $t\in S(Q)$ there is a unique $\varphi\in\cA(Q)$ with
$\varphi(c_Q)=t$.
\item{(2)}  If $\varphi_1,\varphi_2\in\cA(Q)$ satisfy
$\varphi_1(c_Q)\ne\varphi_2(c_Q)$, then $\varphi_1(y)\ne\varphi_2(y)$ for all
$y\in Q$.
\item{(3)} For any $\epsilon>0$, there exists $J$ and $\delta>0$ such that if
$j\ge J$, $\varphi_1,\varphi_2\in\cA(Q,n_j)\cup\cA(Q)$ satisfy
$|\varphi_1(c_Q)-\varphi_2(c_Q)|<\delta$, then
$||\varphi_1-\varphi_2||_Q<\epsilon$.

\give Proof.  We will prove (1); the assertion (2) follows from the Hurwitz
Theorem, and (3) then follow from (2.6).  Let us suppose that there are
distinct functions
$\varphi_1$ and
$\varphi_2\in\cA(Q)$ with $\varphi_1(c_Q)=\varphi_2(c_Q)$.  By Lemma 6.4 of
[BLS1], we may move the point $y=c_Q$, if necessary, to have
$\varphi_1'(c_Q)\ne\varphi_2'(c_Q)$.  Let us write $t^{(k)}=\varphi_k'(c_Q)$,
for $k=1,2$ and set $\epsilon=|t^{(1)}-t^{(2)}|$.  

Let $\{n^{(k)}_j\}$  denote the subsequence of $\{n_j\}$ which produced
$\varphi_k$.  Now by (2.6) it follows that there is a neighborhood
$U$ of $(\varphi_1(c_Q),c_Q)$ and a large number $J$ such that if $j\ge J$,
then for any graph $\Gamma$ from $\cA(Q,n_j^{(k)})$, the slope of
$\Gamma$ is within $\epsilon/2$ of $t^{(k)}$ at all points of $U\cap \Gamma$. 
But this is a contradiction, for if we write $\mu^-_Q$ in polar form, as a
tangent 2-vector times a measure, then on $U$ the tangent vector must be within
$\epsilon/2$ of both $t^{(1)}$ and $t^{(2)}$.  \qed

Passing to further subsequences, we may assume that $\cA(Q_1)\subset\cA(Q_2)$
if $Q_1\in\cG_{s_1}$, $Q_2\in\cG_{s_2}$, and $Q_1\supset Q_2$.  Thus if we write
$\mu^-_{\cQ_s}=\sum_{Q\in\cQ_s}\mu^-_Q,$
then $\mu^-_{\cQ_s}\le\mu^-_{\cQ_{s+1}}$.

\proclaim Theorem 2.4.  The currents $\mu^-_{\cQ_s}$ increase to $\mu^-$ as
$s\to\infty$.  Further each $\mu^-_Q$ has a uniform laminar structure given
by $\mu^-_Q=\int_{a\in S(Q)}[\Gamma_a]\,\nu_Q(a)$

\give Proof.  We have already that
$\mu^-_{Q,n}=\int_{a\in S(Q,n)}[\Gamma_a]\,\nu_Q(y,n)$ for any $y\in Q$.  Now
as $j\to\infty$ we have $\nu_Q(y,n_j)\to\nu_Q$ by Lemma 2.2 and
$\cA(Q,n_j)\to\cA(Q)$ by Lemma 2.3, and thus the integral representations
converge.  This proves that $\mu^-_Q$ has the uniform laminar structure.  
We know that $k^{-1}d^{n_0-n}[f^nX]$ converges to $\mu^-$ and $\mu^-_{Q,n}$
converges to $\mu^-_Q$ as $n\to\infty$.  Thus the inequality $\mu^-_{Q,n}\le
k^{-1}d^{n_0-n}[f^nX]$ yields  $\mu^-_Q\le\mu^-$ and thus
$\mu^-_{\cQ_s}\le\mu^-$.  Similarly, the estimate in Lemma 2.1 converges to
${\bf M}[(\mu^--\mu^-_{\cQ_s})\contract\pi^{-1}Q_0]\le Cd^{-2s}m(\kappa)$.
Thus $\lim_{s\to\infty}\mu^-_{\cQ_s}=\mu^-$. \qed

Our derivation of laminar structure up to this point has relied on the fact
that, as a current, $\mu^-$ has complex dimension 1, and thus sets of area zero
inside each leaf are invisible from the point of view of $\mu^-$.  For the
purpose of defining the critical measure, we will need to know that this laminar
structure actually has leaves which are ``complete,''  i.e.\ conformally
equivalent to $\cx{}$, since the critical points occur on a discrete subset of
the leaf.  This will be done in \S3.

\section 3.  Pesin Theoretic Properties of the Stable/Unstable Currents

We discuss some results from smooth ergodic theory that we will apply to the
structure of the currents $\mu^\pm$ and the measure
$\mu$.  This includes the existence of Lyapunov exponents and the Pesin theory
for stable/unstable manifolds.  These dynamical methods lead us again to 
``laminar" properties of the stable/unstable currents $\mu^\pm$ with respect to
the Pesin stable/unstable manifolds.  In Lemma 3.3 we show that almost every
leaf in the laminar structure obtained in Theorem 2.4 already contains the Pesin
unstable manifolds.

We define Lyapunov exponents for an ergodic measure $\mu$ for a
diffeomorphism in dimension 2.  By the Oseledets Theorem, there is a measurable,
$f$-invariant complex splitting $E^s_x\oplus E^u_x$ of the tangent space for
$\mu$ almost every point $x$, and there exist numbers $\lambda^s\le\lambda^u$,
such that the limits
$$\lambda^s=\lim_{k\to\pm\infty}{1\over k}\log|Df^k|_{E^s}|,\qquad
\lambda^u=\lim_{k\to\pm\infty}{1\over k}\log|Df^k|_{E^u}|\eqno(3.1)$$
exist.  In particular,  the matrix norm satisfies
$$\lim_{k\to+\infty}{1\over k}\log||Df^k_x||=\lambda^u$$ 
for almost every $x$.

In [BS3] we showed that the Lyapunov exponents of the invariant measure $\mu$
satisfy $\lambda^s\le-\log d<0<\log d\le \lambda^u$; since these are nonzero,
$\mu$ is a {\it hyperbolic measure}.  In the following, we may assume more
generally that $\mu$ is a hyperbolic measure of saddle type: i.e.\ the Lyapunov
exponents satisfy $\lambda^s<0<\lambda^u$. 

Let us recall the set $\cR$ of Oseledets regular points.  General
references for the Oseledets Theorem and the Pesin Theory are Pugh and Shub [PS]
and Pollicott [P]. A point $x$ belongs to $\cR$ if for each $\epsilon>0$, there
is a constant $\gamma_{x,\epsilon}>0$ such that
$$|\lambda^{k,s}_x|=\left|Df^k|_{E^s_x}\right|\le
\gamma_{x,\epsilon}e^{k(\lambda_s+\epsilon)}\eqno(3.2)$$
$$|\lambda^{k,u}_x|=\left|Df^{-k}|_{E^u_x}\right|\le
\gamma_{x,\epsilon}e^{-k(\lambda_u-\epsilon)}\eqno(3.3)$$ $${\rm\
angle}(E^s_{f^kx},E_{f^kx}^u)\ge
\gamma_{x,\epsilon}^{-1}e^{-|k|\epsilon}.\eqno(3.4)$$
By the Oseledets Theorem, $\cR$ is a Borel set of full $\mu$ measure.
This means that we have strict contraction in the inequalities (3.2)
and (3.3) if $\epsilon$ is small.

We note that the mapping $f$ is said to be {\it uniformly hyperbolic} if
inequalities (3.2) and (3.3) hold for some uniform constants $\gamma
e^{-kc}$,
independent of the point $x\in J$.  It follows in the uniform case that
the angle is bounded below, independently of $k$ and $x$. In general, uniformly
hyperbolic diffeomorphisms are quite well behaved.

A result of the Pesin theory is that for each regular point $x\in\cR$
the set
$$\eqalign{
W^s(x)&=\{q\in\cx2:\lim_{n\to\infty}\hbox{\rm dist}(f^nq,f^nx)=0\}\cr
&=
\{q\in\cx2:\lim_{n\to\infty}{1\over n}\log\hbox{\rm dist}
(f^nq,f^nx)=\lambda^s\}}$$
is a 2-dimensional  imbedded submanifold.  In the complex case, $W^s(x)$
is a complex manifold (Riemann surface). For $\mu$ almost every $x\in\cR$
the manifold $W^s(x)$ is conformally equivalent to $\cx{}$ (see [BLS1,
Proposition 2.6] or [W]).

Let us consider a coordinate chart $\psi:U\to\Delta^2=\{|x|,|y|<1\}$ for some
open set $U\subset\cx2$, and let us work on $\Delta^2$.  An analytic graph
$T=\{x=g(y):y\in\Delta\}$ will be called a vertical transversal; we
define a horizontal transversal similarly.  We will define a {\it stable
box} $B^s$ (with respect to $\Delta^2$) to be a union of components
$\Gamma$
of $W^s(x)\cap\Delta^2$ for $x\in\cR$ such that $\Gamma$ is a horizontal
transversal to $\Delta^2$.  Thus for any vertical transversal
$T\subset\Delta^2$
there is a set $E\subset T$ such that $B^s=\bigcup_{t\in E}\Gamma_t^s$,
where $\Gamma_t^s$ is a horizontal transversal such that
$\Gamma_t^s\subset
W^s(x)$ for some $x\in\cR$, and the point $t$ is defined by
$\{t\}=\Gamma^s_t\cap T$.  It follows that distinct $\Gamma_t^s$ are pairwise
disjoint, and
$t\mapsto\Gamma_t^s$ is continuous.  We define an unstable box
$B^u=\bigcup_{t\in E^u}\Gamma_e^u$ in a similar fashion, with the unstable
disks taken to be
vertical transversals.

If $B^s$ and $B^u$ are stable and unstable boxes in the same coordinate
neighborhood $\Delta^2$, then the intersection $B=B^s\cap B^u$ is called
a {\it Pesin box}.  The stable and unstable manifolds give $B$ the
structure
of a topological product.  By
[BLS1, Theorem 4.7] the restriction
$\mu\contract B$ is the product measure $\tau^s\otimes\tau^u$ with
respect to this
topological product structure. 

 By the Pesin theory, $\cR$ may be
covered, up to a
set of $\mu$ measure zero, by a countable family of Pesin boxes
$\{B_j\}$.
Thus for a.e.\ $x$ there exists $\epsilon(x)>0$ such that
$W^s_{loc}(x,\epsilon(x))$ is contained in the unstable box $B^u_j$
associated with the Pesin box $B_j$.  We let
$\tilde\cR$ denote the
points $x\in\cR$ such that
$$W^s(x)\subset\bigcup_{n\le 0}f^n\Big(\bigcup_jB_j^u\Big).$$ The following
allows us to ignore the subset of  $\cW^s$ which is not covered by stable
boxes.

\proclaim Proposition 3.1.  $\tilde \cR$ is an $f$-invariant set of full $\mu$
measure.
Further, $$\bigcup_{x\in\cR}W^s(x)-\bigcup_{x\in\tilde\cR}W^s(x)$$ has
zero measure for every slice measure $\mu^+|_T$.  (And thus this
set has zero $|\mu^+|$-measure.)

\give Proof. Almost every point $x$ is contained in a stable box $B^s$, and
there is a number $r(x)$ such that the  stable leaf in $B^s$ containing $x$
is the graph over a Euclidean disk of radius $r(x)$ and centered at $x$.  For
$0<\epsilon$ and $C<\infty$, we let $S(\epsilon,C)$ denote the set of points
$x$ such that $r(x)\ge\epsilon$ and (3.2--4) holds for $\gamma_{x,\epsilon}\le
C$.  By choosing
$C$ large and $\epsilon$ small, we have $\mu(S(\epsilon,C))>0$.  By Poincar\'e
recurrence, almost every $x$ has the property that $f^{n_j}x\in
S(\epsilon,C)$ for infinitely many $n_j\to\infty$.  Let $x$ be such a point and
set $x_j=f^{n_j}(x)$.  Without loss of generality, we may assume that
$\epsilon=1$.  Let $D_j$ denote a copy of the unit disk, and let
$\chi_j:D_j\to W^u(x_j)$ be a conformal coordinate chart with
$\chi_j(0)=x_j$ which expresses the local stable manifold as a graph over
$D_j$ in coordinates such that the graph is flat to first order over the origin.

Now we have a family of germs of conformal mappings
$\varphi_j:D_j\to D_{j+1}$ of a neighborhood of the origin which satisfy
$\varphi_j=\chi_{j+1}^{-1}\circ f^{n_{j+1}-n_j}\circ\chi_{j}$.  Thus
$\varphi_j(0)=0$ and
$|\varphi'(0)|=\left|Df^{n_{j+1}-n_{j}}|_{E^s(x_j))}\right|$.

Given $0<\rho<1$, we choose $\kappa$ such that $\kappa<(1-\rho)^2/2$.  We may
pass to a subsequence of $\{n_j\}$ so that $n_{j+1}-n_j\to\infty$ arbitrarily
fast. By (3.3) we may assume that $|\varphi'_j(0)|\le\kappa$ for each
$j$.  We let $D_{j,\rho}$ denote the disk of radius $\rho<1$ inside $D_j$.
For $R>0$ sufficiently small, $\varphi_j$ is defined on $D_{j,R}$, and by
the Distortion Theorem in one complex variable, the image $\varphi_j(D_{j,\rho
R})$ is contained in the disk of radius $\rho R|\varphi_j'(0)|/(1-\rho)^2$.  By
the choice of
$\kappa$, it follows that $\varphi_j$ extends to all of $D_{j,\rho}$, and
$\varphi_j(D_{j,\rho})$ is contained in the disk of radius
$\rho|\varphi_j'(0)|/(1-\rho)^2$.
This number is less than
$\rho/2$, and so the modulus of the annulus
$D_{j+1,\rho}$ minus the closure of $\varphi_j(D_{j,\rho})$ is at least
$\log 2$.

Now we define
$$W:=\bigcup_{j=1}^\infty f^{-n_j}D_j\subset W^u(x).$$
It follows that $W$ is the increasing union of annuli of moduli at least $\log
2$, and so $W$ is conformally equivalent to $\cx{}$.  Thus $W=W^s(x)$.  It
follows that $x\in\tilde\cR$, and so the $\mu$ measure of $\cR-\tilde\cR$ is
zero.

The statement concerning the slice measures $\mu^+|_T$ follows because $\mu$
has a local product structure, with the factors given by the stable slices
$\mu^+|_{T'}$ and the unstable slices $\mu^-|_{T''}$.
\qed

If $B^s$ is a stable box in the bidisk $\Delta^2$, and if $T$ is a
vertical transversal to $\Delta^2$, then the restriction
$\mu^+|_T\contract(T\cap B^s)$
of the induced measure to $T\cap B^s$ will be called a transversal measure.
For two vertical transversals $T_1$ and $T_2$ of $\Delta^2$, there is a
homeomorphism
$\chi_{}:T_1\cap B^s\to T_2\cap B^s$
obtained by following an intersection point $t_1=T_1\cap \Gamma_t$ along
the graph of a stratum $\Gamma_t$ to the intersection point $T_2\cap\Gamma_t$.
By [BLS1, Theorem 4.5], $\chi_{}$ preserves the set of transversal measures:
$$(\chi_{})_*(\mu^+|_{T_1}\contract(B^s\cap T_1)) =
\mu^+|_{T_2}\contract(B\cap T_2).\eqno(3.6)$$
If $B^s=\{\Gamma_t:t\in E\}$ is a stable box, then in [BLS1] 
the restriction of
$\mu^+$ to $B^s$ was shown to be equal to
$$\mu^+\contract B^s=\int\mu^+_\tau(t)[\Gamma^s_t]\eqno(3.7)$$
where $\mu^+_\tau$ is any transversal measure.   Likewise, for an unstable box
$B^u$, we have a similar representation for
$\mu^-\contract B^u$.
The transformation rule $f_*\mu^+=d^{-1}\mu^+$ corresponds to the fact
that the push-forward under $f_*$ of a transversal measure is $1/d$ times
another transversal measure.  

We may define a wedge product $dd^cU\wedge
T$ for any bounded, continuous psh function $U$ and positive, closed
current $T$, where if $\xi$ is any test form, the product $U\xi$ is a
compactly supported form with continuous coefficients, so we may set
$$dd^cU\wedge T(\xi):= T(Udd^c\xi)$$ (see [BT] for further discussion of this
wedge operation on currents).  A related operation is the intersection
product, $[Z_1]\wedge[Z_2]$, which gives the current of integration over the
intersection $[Z_1\cap Z_2]$.  By integration with respect to the transversal
measure, we may define an intersection wedge product $\dot\wedge$ of a current
of the form (3.7) and a current of integration $[Z]$.  In [BLS1, Lemma 8.3] it
was shown that if
$Z$ is a complex variety, and if $\mu^+$ has the form (3.7) then these two
notions of wedge product coincide, i.e. 
$$(\mu^+\contract B)\wedge[Z]=\int\mu^+_\tau(t)\,[\Gamma_t\cap Z].\eqno(3.8)$$
Because of this, we will use intersection products whenever it is convenient,
but we will just use the notation $\wedge$.

\proclaim Proposition 3.2.  There are countably many unstable boxes $B_j^u$
such that the splitting $E^s\oplus E^u$ extends continuously to $B_j^u$, there
is a constant $C_j$ such that (3.3--4) hold on $B_j^u$ with
$\gamma_{x,\epsilon}\le C_j$ for $x\in B_j^u$, and such that for any complex
manifold $T$,
$T\cap\bigcup_{j,n=1}^\infty f^nB_j^u$ has full measure for the slice measure
$\mu^-_T$.

\give Proof.  Let $\{B_j^u\}$ be a family of unstable boxes as in
Proposition 3.1.  We may choose stable boxes $B^s_j$ such that $B_j=B^s_j\cap
B^u_j$ is a Pesin box, and $E^{s/u}$ extend continuously to 
$B_j^u$ and (3.3--4) hold.   \qed

Another consequence of the Pesin theory is that there is a measurable family of
Lyapunov charts.  This means that almost every $x$ is the center of a
(complex) affine image $L(x)$ of a bidisk $\Delta^2$, and there is a product
metric on $L(x)$ which is strictly expanded/contracted under $f$ (see [PS]). 
(We call $L(x)$ a topological bidisk in [BLS2].)  If $X$ is a complex variety,
the cutoff image of $X$ under $f$, i.e.\ $f(X\cap L(x))\cap L(fx)$ is
stretched across $L(fx)$ in the unstable direction.  In fact, if $X\cap L(x)$
intersects $W^s_{loc}(x)$ transversally at $x$, then there is a number $N(x)$
such that if $m\ge N(x)$, then after $m$ stretchings and cuttings-off, we have an
unstable transversal to
$L(f^mx)$, i.e.\ 
$$f^m(X\cap L(x))\cap f^{m-1}L(fx)\cap\cdots\cap L(f^mx)\eqno(3.9)$$
is an unstable transversal to $L(f^mx)$.

Let us take a countable family of Pesin boxes $B_j$ whose union has full
measure and which have the property that the constant $\gamma_{x,\epsilon}$ in
(3.2--4) satisfies $\gamma_{x,\epsilon}\le C_j$ for $x\in B_j$.  Further, we may
assume that the inner radius of $L(x)$ is bounded below by $r_0>0$ for all $x\in
B_j$.  Further, we may assume that the axes of the bidisk $L(x)$ are almost
constant for $x\in B_j$, and we may assume that the projection $\pi$ is
transversal to the unstable direction, i.e.\ $\pi^{-1}(0)$ makes a positive
angle with the unstable axis of $L(x)$.  Srinking $B_j$ if necessary, we may
assume that there is a square $Q_j$ with $B_j\subset\pi^{-1}Q_j$, and such
that $\pi^{-1}q\cap L(x)$ is a vertical transversal of $L(x)$ for all $x\in B_j$
and $q\in\bar Q_j$.  Finally, $W^u_{loc}(x)\cap L(x)$ is an unstable transversal
to $L(x)$, so we may assume that for each stratum $\Gamma$ of $B_j^u$, $\Gamma$
crosses $\pi^{-1}Q_j$ properly, i.e.\ the restriction of $\pi$ from
$\Gamma\cap\pi^{-1}Q_j$ to $Q_j$ is a homeomorphism.

\proclaim Lemma 3.3.  There are countably many Pesin boxes $\{B_j\}$ such
that $\bigcup B_j$ has full $\mu$ measure, and for each $B_j$ there is a
square $Q_j\subset\cx{}$ such that the associated unstable box $B^u_j$
satisfies 
$$\mu^-_{Q_j}\ge\mu^-\contract (B^u_j\cap\pi^{-1}Q_j).$$

\give Proof.  We take $B_j$ and $Q_j$  as in the discussion above.  We note
that we may take $B_j^u$ such that for each stratum $\Gamma$ of $B^u_j$,
$\Gamma\cap L(x)$ is an unstable transversal to $L(x)$ for all
$x\in B_j$. Let $\{P_j\}$ be a finite
family of disjoint Pesin boxes with a
family of disjoint open sets $V_j$ with $V_j\supset P_j$.  Further, if
$\epsilon>0$ is given, we may assume that $\mu(\bigcup P_j)>1-\epsilon$. 

For a fixed $j$, we will set $B=B_j$ and $Q=Q_j$ and will show that they
have the property claimed in the Lemma.  Let
$c>0$ be such that $\lim_{n\to\infty}cd^{-n}[f^nX]=\mu^-$.  We may suppose that
$$\int(\mu^+\contract P^s_j)\wedge c[X\cap V_j]\ge
(1-\epsilon)\mu(P_j),\eqno(3.10)$$
replacing $X$ by $f^nX$ and $c$ by $cd^{-n}$, $n$ large, if necessary.  For
each $m$, let $\cG(m,j)$ denote the set of connected components $\Gamma$ of
$f^m(X\cap V_j)\cap\pi^{-1}Q$ such that $\pi|_\Gamma:\Gamma\to Q$ is a
homeomorpism.  We let
$$\mu^-_{\cG(m,j)}=cd^{-m}\sum_{\Gamma\in\cG(m,j)}[\Gamma]\eqno(3.11)$$
so that
$$\mu^-_Q\ge\limsup_{m\to\infty}\sum_j\mu^-_{\cG(m,j)}.$$
The inequality arises since there are possibly good disks in $f^m(X)\cap
\pi^{-1}Q$ that are lost when $f^m(X-\bigcup V_j)$ is removed.  We note that
since each $\Gamma\in B_j^u$ is a proper transversal to $\pi^{-1}Q$, it will
suffice to show that
$$\mu^-_Q\wedge(\mu^+\contract B^s)\ge(\mu^-\contract
B^u)\wedge(\mu^+\contract B^s)=\mu(B).\eqno(3.12)$$
In (3.12) it is the inequality that needs to be proved; the equality is just
the product structure of $\mu$ on $B$.

Now we have
$$\eqalign{
\mu^-_Q\wedge (\mu^+\contract
B^s)&\ge\limsup_{m\to\infty}\sum_j\mu^-_{\cG(m,j)}\wedge(\mu^+\contract
B^s)\cr 
&\ge\sum_j(1-\epsilon)\mu(P_j)\mu(B)\cr
&\ge(1-\epsilon)^2\mu(B),\cr}$$
where the second inequality follows from Lemma 3.4 below.   Thus we conclude
that (3.12) holds, which completes the proof. \qed

\proclaim Lemma 3.4.  Let $B$ and $P_j$ be as above. Then
$$\lim_{m\to\infty}\int\mu^-_{\cG(m,j)}\wedge(\mu^+\contract B^s)
\ge(1-\epsilon)\mu(P_j)\mu(B).$$

\give Proof.  We note that each unstable transversal $\Gamma$ in $L(x)$ gives
rise to a unique good disk $\Gamma\cap\pi^{-1}Q$.  Thus we will consider
instead the current $\mu^-_{\cV(m,j)}$, where the sum in (3.11) is replaced by
$\Gamma\cap\pi^{-1}Q$ for $\Gamma$ which are unstable transversal components
of $f^m(X\cap V_j)\cap L(x)$ for some $x\in B$.  Since
$\mu^-_{\cG(m,j)}\ge\mu^-_{\cV(m,j)}$, it suffices to prove the Lemma for
$\mu^-_{\cG(m,j)}$ replaced by $\mu^-_{\cV(m,j)}$.  By Lemma 6.4 of [BLS1], we
may suppose that $X$ intersects $W^s_{loc}(x)$ transversally for each $x\in
P_j$.  Let $N(x)$ denote the measurable function on $X\cap V_j$ with the
property (3.9).

We define 
$c_1=\int\mu^+\wedge c[X\cap V_j]$
so that  $c_1\ge(1-\epsilon)\mu(P_j)$ by (3.10).  We may assume (changing
$V_j$ slightly if necessary) that $\mu^+|_X$ puts no mass on $\partial(X\cap
V_j)$.  Thus
$$\lim_{m\to\infty} cd^{-m}[f^m(X\cap V_j)]=c_1\mu^-.$$
It follows that
$$\eqalign{
\lim_{m\to\infty}\int cd^{-m}[f^m(X\cap V_j)]&\wedge(\mu^+\contract
B^s) =  c_1\int\mu^-\wedge(\mu^+\contract B^s)\cr
&=c_1\mu(B)\ge(1-\epsilon)\mu(P_j)\mu(B).\cr}$$
Thus if we set
$$\eta^-(m,j)=cd^{-m}[f^m(x\cap V_j)]-\mu^-_{\cV(m,j)}$$
it will suffice to show that
$$\lim_{m\to\infty}\int\eta^-(m,j)\wedge(\mu^+\contract B^s)=0.$$
However, if we pull back to $X\cap V_j$ via $f^m$ and recall the definiton of
$N(x)$, we have
$$\int\eta^-(m,j)\wedge(\mu^+\contract B^s)=\int_{\{N(x)>m\}}[X\cap
V_j]\wedge(\mu^+\contract P^s_j).$$
Thus the right hand side tends to 0 as $m\to\infty$ since
$\{N(x)>m\}$ decreases to $\emptyset$. \qed

\section 4.  Averaged Rates of Growth

    Lyapunov exponents describe the behavior tangent vectors at
$\mu$ a.e.\ point. This is not however the most direct way to get hold of the
value of the Lyapunov exponents.  In this Section we consider various
alternative notions of the growth rate of vectors  and we relate them to the
Lyapunov exponent.  We  discuss a method of measuring the growth of
$Df^k$ by taking the average with respect to $\mu$ and all directions; and
we show how it is related to a type of critical point.  Finally, we give a
formula for the averaged rate of growth by pulling back a form from
projective space.  This last formula (Proposition 4.6) is of interest because
it involves the projectivized image of the map
$x\mapsto Df^n_x\in\cL(\cx2,\cx2)$, and thus measures the volume of the
(projectivized) image rather than the size of $\Vert Df^n\Vert$. This
description suggests an analogy with the definition of curvature via
the Gauss map.

There is a certain symmetry between $\lambda^+$ and $\lambda^-$ which can be
realized by replacing $f$ by $f^{-1}$. For the sake of definiteness we focus on
$\lambda^+$ in this Section, and our notation reflects that emphasis. 

 We let
$\alpha,\beta$ denote constant, nonzero vector fields on
$\cx2$, and define the quantities
$$\Lambda=\lim_{k\to\infty}{1\over k}\int\log||D^kf(x)||\,\mu(x).$$
$$\Lambda(\alpha)=\lim_{k\to\infty}{1\over
k}\int\log|Df^k_x(\alpha)|\mu(x)$$
$$\Lambda(\alpha,\beta)=
\lim_{k\to\infty}{1\over k}\int\log|\beta\cdot Df^k_x(\alpha)|\mu(x).$$

The first integral arises in the proof of the Oseledec Theorem as the first
step in the proof of the existence of Lyapunov exponents. From this we see that
$\Lambda=\lambda^+(\mu)$. We will analyze the other two quantities.

We can identify $\alpha$ and $\beta$ with vectors in $\C^2-\{0\}$. It is clear
that $\Lambda(c\alpha)=\Lambda(\alpha)$ for any $c\in\cx{}-\{0\}$ and
$\Lambda(c_1\alpha,c_2\beta)=\Lambda(\alpha,\beta)$ for $c_1,c_2\in\cx{}-\{0\}$.
Thus when it is convenient we may think of identify $\Lambda(\cdot)$ and
$\Lambda(\cdot,\cdot)$ as functions on $\cp1$ and $\cp1\times\cp1$.  On the other
hand it is sometimes convenient to assume that $|\alpha|=1$ and $|\beta|=1$.  We
let $\sigma$ denote the rotation invariant measure on the unit ball in $\cx2$,
normalized to have total mass 1. We denote by the same letter the induced
measure on $\cp1$.

In the sequel, we will use the observation that
if $\beta=(\beta_1,\beta_2)\in\cx2$, then
$$\int_{|\alpha|=1}\log|\alpha\cdot
\beta|\,\sigma(\alpha)=\log|\beta|-{1\over2}\eqno(4.1)$$
depends only on $|\beta|$.

\proclaim Lemma 4.1.  For $\sigma$ a.e.\ $\beta\in\cp1$,
$\Lambda(\alpha,\beta)=\Lambda(\alpha)$.

\give  Proof.  For each $x$ we have
$$\int_{\beta\in\cp1}\log|\beta\cdot
Df^k_x(\alpha)|\sigma(\beta) =\log|Df^k_x(\alpha)|-{1\over2}$$
by (4.1).  Now we integrate with respect to $\mu(x)$, divide by $k$, and then
take the limit as $k\to\infty$ to obtain
$$\int_{\beta\in\cp1}\Lambda(\alpha,\beta)\sigma(\beta)=\Lambda(\alpha).$$

On the other hand, we may assume that $|\beta|=1$, so 
$\Lambda(\alpha)\ge\Lambda(\alpha,\beta)$.  Thus the Lemma follows.\qed

Recall from \S3 the measurable, $f$-invariant complex splitting $E^s_x\oplus
E^u_x$ of the tangent space for $\mu$ almost every point $x$. We let $x\mapsto
e^{s/u}_x$ be a measurable choice of unit vectors in $E^{s/u}_x$.  Given a
tangent vector $\alpha=\alpha_1\partial_1+\alpha_2\partial_2$ (using
$\partial_1$ and $\partial_2$ to denote a frame for the tangent space of
$\cx2$ the point
$x$) we may split it as 
$$\alpha=\alpha_x^se_x^s +\alpha_x^ue_x^u.$$ 
Thus for
$\mu$ a.e.\ $x$ there are numbers $\lambda_x^{k,s/u}$ such that
$$Df^k_{x}(\alpha) = \lambda_x^{k,s}\alpha_x^se_{f^kx}^s +
\lambda_x^{k,u}\alpha_x^ue_{f^kx}^u.$$
Thus we have represented $Df^k$ as a diagonal matrix.

\proclaim Lemma 4.2.  For $\sigma$ a.e.\ $\alpha$, we have
$\Lambda(\alpha)=\Lambda$.

\give Proof.  The function $\alpha_x^s$ in the splitting above is given
by the
Hermitian inner product $\langle\alpha,e_x^s\rangle$ on $\cx2$.  For $x$
fixed,
$\int\log|\langle\alpha,e_x^s\rangle|\sigma(\alpha)=-{1\over2}$ as above.  Since
the integrand is  nonpositive,
it follows that $\log|\alpha^s_x|$ is integrable with respect to the
product
measure $\sigma\times\mu$.  Reversing the order of integration, we have
$$\int\mu(x)\int\log|\langle\alpha,e_x^s\rangle|\sigma(\alpha) =
\int\int\mu(x)\log|\langle\alpha,e_x^s\rangle|\sigma(\alpha).$$
Thus for almost
every $\alpha\in\cp1$ the function $x\mapsto\log|\alpha_x^s|$ is
integrable with
respect to $\mu$. Similarly, we may assume that $\log|\alpha_x^u|$ is
integrable.

Letting $\gamma_{x,\epsilon}$ be as in (3.2--4), we define
$S_\gamma=\{x:\gamma_{x,\epsilon}\le\gamma\}$ for fixed $\epsilon>0$.
By the splitting above, we have
$$\log|Df^k_x(\alpha)| = \log|\lambda_x^{k,s}\alpha_x^se_{f^kx}^s +
\lambda_x^{k,u}\alpha_x^ue_{f^kx}^u| = \log|A+B|.$$
Given two vectors $A$ and $B$ which form an angle of opening $\theta$, the
square of the sum has length
$$\eqalign{
|A+B|=&|A|^2+|B|^2+2|A|\,|B|\cos\theta\cr
=&(|A|-|B|)^2+2(\cos\theta-1)|A|\,|B|\ge 2(\cos\theta-1)|A|\,|B|.\cr}$$
For $\theta$ small, we may estimate $2(\cos\theta-1)$ by $\theta^2$, so
we have
$$\log|Df^k(\alpha)|
\ge \log\left|\theta^{2}\alpha_x^s\lambda_x^{k,s}
\alpha_x^u\lambda_x^{k,u}\right|$$
Thus by (3.2),
(3.3), and (4.1) and the fact that the angle between $e^s_{f^kx}$ and
$e^u_{f^kx}$ is bounded below by (3.4) and $\gamma\ge\gamma_{x,\epsilon}$
for
$x\in S_\gamma$ we have that the quantity
$$\eqalign{ {1\over k}\log|Df^k(\alpha)|
\ge & \log\theta^2 + \log|\lambda_x^{k,s}\lambda_x^{k,u}| +
\log|\alpha_x^s| + \log|\alpha_x^u|\cr
\ge & (-2\epsilon-{2\over k}\log\gamma) +(\log|a|-2\epsilon -
{2\over k}\log\gamma) + \log|\alpha_x^s| + \log|\alpha_x^u| \cr}$$
is bounded below by a function which is integrable with
respect to $\mu$.

For $\mu$ a.e.\ point $x$ such that $\alpha_x^u\ne0$, we have
$$\lim_{k\to\infty}{1\over k}\log|Df^k_x(\alpha)|=\Lambda,$$
so by (4.1) and the dominated convergence theorem, we have
$$\lim_{k\to\infty}\left|\int(\Lambda -{1\over
k}\log|Df^k_x(\alpha)|)\mu(x)\right|\le \mu(J-S_\gamma)(\Lambda+\log
M).$$
The Lemma follows since $\lim_{\gamma\to\infty}\mu(J-S_\gamma)=0$. \qed

\proclaim Lemma 4.3.  For a.e.\ $\beta$ we have
$$\lim_{k\to\infty}{1\over k}\int\sigma(\alpha)\int\log|\beta\cdot
Df_x^k(\alpha)|\mu(x)=\Lambda.$$

\give Proof.  This follows from Lemma 4.1, 4.2 and the bounded
convergence theorem.  \qed

Another way in which a family of critical points arises is as follows.
Let us define
$$Z_k(\alpha,\beta)=\{x\in\cx2:\beta\cdot Df^k_x(\alpha)=0\}\eqno(4.2)$$
as the critical points of the scalar function $x\mapsto\beta\cdot f$
with respect to the direction $\partial_\alpha$.  Unlike the set of unstable
critical points $\cC^u$, which will be defined in \S5, this set is is not
invariant.  On the other hand, it is quite explicit.

The following  computation resembles the proof of [BS3, Theorem 3.2],
except that now we keep the integral term on the right hand side for further
analysis.

\proclaim Lemma 4.4.  Let $\alpha,\beta\in\cx2$ be such that the second
coordinates $\alpha_2,\beta_2\ne0$ are nonvanishing.  If $T=\{x=0\}$ is the
$y$-axis, then
$$\int\log\left|\beta\cdot Df^k\alpha\right|\,\mu^+\wedge {1\over{d^n}}
f^n_*[T] =
\log|\alpha_2\beta_2 d^k| +\int G^+[Z_k(\alpha,\beta)]\wedge
{1\over{d^n}}f^n_*[T].$$

\give Proof.  Applying $(f^n)^*$ and treating $f^{n*}f^n_*$ as the identity
transformation, we have
$$
\int\log|\beta\cdot Df^k\alpha|{1\over 2\pi}dd^cG^+\wedge {1\over{d^n}}f^n_*[T]
=
\int f^{n*}\left(\log|\beta\cdot
Df^k\alpha|\right){1\over2\pi}dd^cG^+\wedge[T],$$
where we use the functional equation $f^{n*}G^+=G^+\circ f^n=d^nG^+$.
Furthermore,
$G^+$ restricted to $T$ is the Green function of $K^+\cap T$, so that
${1\over2\pi}(dd^c)_TG^+|_T$ is the harmonic measure, which we denote by
$\mu^+_{K^+\cap T}$, so the equation becomes
$$\int\log|\beta\cdot Df^k\alpha|{1\over 2\pi}dd^cG^+\wedge
{1\over{d^n}}f^n_*[T] =\int_Tf^{n*}\left(\log|\beta\cdot
Df^k\alpha|\right)\mu_{K^+\cap T}.$$
>From formula (1.2) for $f^k$, we  observe that
$$\beta\cdot Df^k(\alpha)=\beta_2\alpha_2d^ky^{d^k-1}+\dots\ .$$
Since $\mu_{K^+\cap T}$ is harmonic measure, we may apply Jensen's
formula [BS3, Lemma 3.1] to the monic polynomial
$(\beta_2\alpha_2d^k)^{-1}\beta\cdot Df^k(\alpha)$ (restricted to $T$) and obtain
$$\eqalign{
= & \log|\beta_2\alpha_2d^k| +\sum_{\{c\in T:\beta\cdot
Df^k_{f^nc}(\alpha)=0\}}d^nG^+(c)\cr
= & \log|\beta_2\alpha_2d^k| + \int G^+[Z_k(\alpha,\beta)]
\wedge {1\over{d^n}}f^n_*[T],\cr }$$
where the last equation comes from pushing $[T]$ forward under $f^n$.
This gives the desired formula.  \qed

\proclaim Corollary 4.5.  If $\alpha_2$, $\beta_2\ne0$, then
$$\Lambda(\alpha,\beta)=\log d+\lim_{k\to\infty}{1\over k}
\int G^+\mu^-\wedge [Z_k(\alpha,\beta)].$$

\give Proof.  We take the formula given in Lemma 4.4 and let
$n\to\infty$.  Then we divide by $k$ and take the limit as $k\to\infty$.
\qed

Now we find another way to replace the explicit dependence on $\alpha$
and $\beta$ by the average over all directions.  This provides an
alternative approach to critical points. 

 The differential induces a mapping
$$\cx2\ni x\mapsto Df_x\in\cL(\cx2,\cx2).$$
We may identify the dual space
$\cL(\cx2,\cx2)^*\cong\cx2\otimes(\cx2)^*$,
where $\alpha\otimes\beta\in\cx2\otimes(\cx2)^*$ induces the functional
$\cL(\cx2,\cx2)\ni Z\mapsto\beta\cdot Z\alpha$.  Let
$$(\alpha\otimes\beta)^\perp=\{Z\in\cL(\cx2,\cx2):\beta\cdot
Z\alpha=0\},$$
and let $[(\alpha\otimes\beta)^\perp]$ denote the current of integration
over $(\alpha\otimes\beta)^\perp$ as a subset of $\cL(\cx2,\cx2)$.  Now
the
function
$V_{\alpha,\beta}(Z)=\log|\beta\cdot Z\alpha|$ satisfies the Poincar\'e-%
Lelong identity
$$[(\alpha\otimes\beta)^\perp]={1\over 2\pi}dd^cV_{\alpha,\beta}.$$
Averaging the function $V_{\alpha,\beta}$ with respect to
$\alpha$ and $\beta$  we have
$$\tilde V(Z):=\int_{\alpha\in{\bf P}^1}\sigma(\alpha)
\int_{\beta\in{\bf P}^1}\sigma(\beta)\, V_{\alpha,\beta}(Z),\eqno(4.3)$$
so that $\tilde V(Z)$ is continuous off the origin, plurisubharmonic, and
logarithmically
homogeneous.  Observe that the integral on the left
hand side, as a function of $Z$, is invariant under the $U(2)\times
U(2)$-action
on $\cL(\cx2,\cx2)$ given by
$(S,T)\cdot Z\mapsto SZT^{-1}$.  Thus to evaluate the integral it suffices to
consider the case where $Z=\hbox{\rm diag}\{\lambda_1,\lambda_2\}$ is diagonal.

In this case $V_{\alpha,\beta}(Z)=
\log|\alpha_1\lambda_1\beta_1+\alpha_2\lambda_2\beta _2|$, and by (4.1) the first
integration inside (4.3) yields
$$\int_\alpha\log\big|\alpha_1\lambda_1\beta_1 +
\alpha_2\lambda_2\beta_2\big|\, \sigma(\alpha)=
{1\over2}\log(|\beta_1\lambda_1|^2+|\beta_2\lambda_2|^2)+C.$$

The (1,1) form
$$\Theta={1\over 2\pi}dd^c\tilde V\eqno(4.4)$$
on $\cx{2\times 2}\cong\cL(\cx2,\cx2)$ represents the averaged
current of integration.  (Note that we are making an abuse of notation,
representing a current as a (1,1)-form.)  By the logarithmic
homogeneity,
we may also interpret  $\Theta$ as a form on the projectivized space
$\cL(\cx2,\cx2)/\cx*\cong
{\bf P}^3$, and $\Theta$ dominates a multiple of the standard K\"ahler
form on
${\bf P}^3$.
If we use again the notation $Df^k$ to denote
the projective image
of the differential in $\cL(\cx2,\cx2)/\cx*$,
then the averaged critical locus is the pullback
of $\Theta$  on projective space:
$$\int_{\alpha\in{\bf P}^1}\int_{\beta\in({\bf
P}^1)^*}\sigma(\alpha)\sigma(\beta)\,[Z_k(\alpha,\beta)]
=(Df^k)^*\Theta.\eqno(4.5)$$

Averaging the formula of Corollary 4.5 over $\alpha,\beta\in{\bf P}^1$,
we obtain:

\proclaim Proposition 4.6.
$$\Lambda=
\log d +\lim_{k\to\infty}{1\over k}\int G^+\mu^-\wedge
(Df^k)^*\Theta.$$

\section 5.  Stable/Unstable Critical Measures

In this Section we begin by defining the unstable critical points $\cC^u$ and the
unstable critical measure $\mu^-_c$. (The definition of the corresponding
objects $\cC^s$ and $\mu^+_c$ should be clear.) We will show  (Theorem
5.1) that if
$\mu^-_{c,s}$ is the critical measure defined starting from the laminar current
$\mu^-_{\cQ_s}$, then $\mu^-_{c,s}$ converges to $\mu^-_c$ as $s\to\infty$.
The rest of the section is devoted to showing (Theorem 5.9) that
$\mu^-_c$ is equal to the limit of the intersection product of
$\mu^-$ with the average over $\alpha$ and $\beta$ of the critical varieties
$f^jZ_k(\alpha,\beta)$ as
$j,k-j\to\infty$, i.e.\ $\mu^-\wedge f^j_*(Df^k)^*\Theta\to\mu^-_c$.

 We define the {\it unstable critical points} as
$$\cC^u=\bigcup_{x\in\cR}\hbox{\rm Crit}(G^+,W^u(x)-K),$$
where $\hbox{\rm Crit}(G^+,W^u(x)-K)$ is the set of critical points, with
multiplicity, of the restriction of the function $G^+$ to the open
subset $W^u(x)-K$ of the manifold $W^u(x)$.  The restriction of $G^+$ to
$W^u(x)$ is subharmonic on $W^u(x)$ and harmonic on $W^u(x)-K$; thus
$W^u(x)-K\ne\emptyset$.  Since $G^+$ vanishes on $W^u(x)\cap K$ (which is
nonempty since it contains $x$), it follows that $G^+$ cannot be constant on a
nonempty open subset of $W^u(x)-K=\{y\in W^u(x):G^+(y)>0\}$.  Further, since
$x\in\cR$ is a regular point, it follows from [BLS1, Proposition 2.9] that the
restriction $G^+|_{W^u(x)}$ is not everywhere harmonic, and so
$W^u(x)-K\ne\emptyset$. Thus $\hbox{\rm Crit}(G^+,W^u(x)-K)$ is a discrete
subset of $W^u(x)-K$ for each $x\in\cR$.
If $f$ is uniformly hyperbolic, then $\cC^u$ is a closed subset of $U^+$. In
the general case $\cC^u$ is likely not to be well behaved.

We will now define the unstable measure $\mu^-_c$. We start by defining its
restriction to an unstable box $B^u$.  For a stratum $\Gamma_t$ of $B^u$, the
critical points of $G^+|_{\Gamma_t-K}$ are discrete, as noted above. We let the
current  $[\hbox{\rm Crit}(G^+,\Gamma_t-K)]$ denote the sum of point masses
(with multiplicity) at the critical points of $G^+|\Gamma_t-K$.  The mapping of
currents $t\mapsto[\hbox{\rm Crit}(G^+,\Gamma_t-K)]$ is semicontinuous and may be
assumed to be bounded, so we may set
$$\mu^-_c\contract B^u=\int\mu^-_\tau(t)\,[\hbox{\rm
Crit}(G^+,\Gamma_t-K)].$$
It is evident that this definition of $\mu^-_c$ is independent of the stable
box involved, since if we have two stable boxes, the two definitions of
$\mu^-_c$ agree on the overlap.  This definition of $\mu^-_c$ may be considered
to give almost all of the points of $\cC^u$, since by Lemma 3.1, we could work
equally naturally with the set $\tilde\cR$, in which case every critical point
would lie inside an unstable box.

Defined this way, $\mu^-_c$ is evidently $\sigma$-finite, and in \S6 we will see
that it is locally finite on $U^+$. The set $\cC^u$ is $f$-invariant. Since the
transversal measures corresponding to $\mu^-$ multiply by $d$ under push-forward
by $f$, and the function $G^+$ multiplies by $d^{-1}$ it follows that
$G^+\mu^-_c$ is $f$-invariant:
$$f_*(G^+\mu^-_c)=G^+\mu^-_c.$$

For a square $Q\in\cQ_s$, we let $\mu^-_Q=\int\nu_Q(a)[\Gamma_a]$ denote the
laminar structure obtained in Theorem 2.4, in terms of an algebraic variety
$X$ and a projection $\pi_\alpha$.  We may define the corresponding critical
measure
$$\mu^-_{c,s}=\sum_{Q\in\cQ_s}\int\nu_Q(a)[{\rm Crit}(G^+,\Gamma_a)].$$

\proclaim Theorem 5.1. For all but countably many values of
$\alpha$,
$$\lim_{s\to\infty}\mu^-_{c,s}=\mu^-_c.$$

\give Proof.  Let us choose Pesin boxes $B_j$ whose union has full $\mu$
measure.   For $\alpha\in\cx2$ we let
$S(\alpha,j)$ denote the set of points of
$B^u_j$ where the tangent space of the corresponding stratum is annihilated by
$\pi_\alpha$.    If $\alpha'$ and $\alpha''$ define different points of ${\bf
P}^1$, then
$S(\alpha',j)$ is disjoint from $S(\alpha'',j)$.  Thus $\mu^-_c(S(\alpha,j))>0$
for only countably many values of $\alpha$.  Thus, except for countably many
values of 
$\alpha$, we have $\mu^-_c(S(\alpha,j))=0$ for all $j$.  Now we may subdivide
the Pesin boxes to obtain a new covering $\{B_j\}$ which staisfies the
hypotheses of Lemma 3.3 with $\pi=\pi_\alpha$.  Thus for each $B_j$ there is a
$Q_j$ such that
$$\mu^-_{Q_j}\ge\mu^-\contract (B^u_j\cap\pi^{-1}Q_j).$$
It follows that if $s$ is sufficiently large that $Q_j$ is a union of squares
from $\cQ_s$, then we have
$$\mu^-_{c,s}\ge\mu^-_c\contract(B^u_j\cap\pi^{-1}Q).$$
It follows  that $\lim_{\kappa\to\infty}\mu^-_{c,\kappa}\contract
Y\ge\mu^-_c\contract Y$ holds for $Y=\bigcup_{j,n}f^n\hat B^u_j$, where $\hat
B^u_j=B^u_j\cap\pi^{-1}Q_j$.  As in Lemma 3.1, $Y$ has full measure with
respect to all transversals, so the Theorem follows. \qed

Now we start the sequence of Lemmas that will lead to the proof of
Theorem 5.9. For $\beta\in(\cx2)^*$, we consider $\alpha\mapsto\beta\cdot
Df^k(z)(\alpha)$ and $\alpha\mapsto\partial G^+\cdot\alpha$ as linear
functionals acting on
$\alpha\in\cx2$.  We let $\langle\beta\cdot Df^k(z)\rangle$ and
$\langle\partial G^+\rangle$ denote their images in $({\bf P}^1)^*$.

\proclaim Lemma 5.2.  For each compact subset $U_0\subset U^+$, the sequence
$\langle\beta\cdot Df^k(z)\rangle$ converges to $\langle\partial G^+\rangle$ as
$k\to\infty$, uniformly in $z\in U_0$ and $\beta\in\cx2-\{0\}$.

\give Proof.   If $f$ has the form (1.1), then the coordinates
$f^n=(f^n_{(1)},f^n_{(2)})$ satisfy
$f^n_{(1)}=(f_2\circ\cdots\circ f_m\circ f^{n-1})_{(2)}$, and
$f^n_{(2)}=(f^n_{(1)})^{d_1}+\cdots$, so $d^{-n}\log|f^n_{(2)}|$ and
$d^{-n}d_1\log|f^n_{(1)}|$ converge to $G^+$ uniformly on compact subsets of
$U^+$.  Thus the normalizations of the gradients $\partial f^n_{(i)}|\partial
f^n_{(i)}|^{-1}$, $i=1,2$ both converge uniformly to the normalization
$\partial G^+|\partial G^+|^{-1}$ on compact subsets of $U^+$.  It follows
that on any compact subset of $U^+$, the normalization of
$\partial(\beta_1f^n_{(1)}+\beta_2f^n_{(2)})$ converges to $\partial
G^+|\partial G^+|^{-1}$ uniformly in $\beta\ne(0,0)$.  Since $\beta\cdot
Df^n(\alpha)$ may be identified with
$\partial(\beta_1f^n_{(1)}+\beta_2f^n_{(2)})\cdot\alpha$, it follows that the
projective images of these linear functionals converge uniformly.\qed

Let $P^u=\bigcup_{t\in T}\Gamma_t$ be an unstable box as in Proposition 3.2. 
For each $j\ge0$, we define 
$$\langle E^s(f^{-j}\Gamma_t)\rangle:=\{\langle E^s_x\rangle :x\in
f^{-j}\Gamma_t\}\subset{\bf P}^1.$$ 
Thus $\langle E^s(f^{-j}\Gamma_t)\rangle$ has diameter $O(e^{-\epsilon j})$. 
We set
$$\cV^s(j,t))=\{\alpha\in{\bf P}^1:\hbox{\rm
dist}(\alpha,\langle E^s(f^{-j}\Gamma_t)\rangle)<{1\over4}\}$$
and $\cV^u(j,t)={\bf P}^1-\cV^s(j,t)$.
It follows from (3.2--4) that $Df^j\cV^s(j,t)$ lies in a
$O(e^{-\epsilon j})$-neighborhood of $\langle E^s_x\rangle$ at all
$x\in\Gamma_t$.

\proclaim Lemma 5.3.  Let $P$ and $\cV^u(j,t)$ be as above, and let us suppose
that $G^+$ has no critical points on $\partial\Gamma_t$ for $t\in T$.  Then
$$\lim_{{j\to\infty}\atop{k-j\to\infty}}{\rm dist}(\Gamma_t\cap
f^jZ_k(\alpha,\beta),{\rm Crit}(G^+,\Gamma_t))=0$$
with the limit being uniform in $t\in T$, $\alpha\in\cV^u(j,t)$, and
$\beta\in({\bf P}^1)^*$.

\give Proof.  If $j,k-j\to\infty$, then there are sequences
$\kappa_1(k),\kappa_2(k)\to\infty$ such that $\kappa_1(k)\le j\le
k-\kappa_2(k)$.    If $\zeta\in\Gamma_t\cap f^jZ_k(\alpha,\beta)$, then
$y=f^{-j}\zeta$ satisfies $\beta\cdot Df^k_y(\alpha)=0$.  Thus $\zeta$
satisfies $\beta\cdot Df^{k-j}_\zeta(f^j_*\alpha)=0$.

For $\delta>0$, we may choose $\kappa_1$ sufficiently large that if
$j\ge\kappa_1$, $\zeta\in\Gamma_t$, and $\alpha\in\cV^u(j,t)$, then
$\hbox{\rm dist}_{{\bf P}^1}(f^j_*\alpha,\langle E^u_\zeta\rangle )<\delta$. 
Furthermore, for $\kappa_2(k)$ sufficiently large and $j\le k-\kappa_2$, it
follows from Lemma 5.2 that $\hbox{\rm dist}_{({\bf P}^1)^*}(\partial
G^+,\beta\cdot Df^{k-j})<\delta$.  Thus the distance between the sets
$\hbox{\rm Crit}(G^+,\Gamma_t)$ and $\{x\in\Gamma_t:\beta\cdot
Df^{k-j}(f^j_*\alpha)=0\}$ is uniformly small. \qed

Next we define
$$\lambda^{s/u}_{j,k}(\beta,t)=\int_{\alpha\in\cV^{s/u}(j,t)
}\sigma(\alpha)[\Gamma_t\cap f^jZ_k(\alpha,\beta)].$$
The plan is to show that $\lambda^u_{j,k}(\beta,t)$ converges to the
critical point measure $[{\rm Crit}(G^+,\Gamma_t)]$, and thus the integral
with respect to $t$ will converge to the critical measure $\mu^-_c\contract
P$, and then to show that $\lambda^s_{j,k}(\beta,t)$ converges to
zero as $j$, $k-j\to\infty$.

\proclaim Lemma 5.4. For each $t\in T$, $\lambda^u_{j,k}(\beta,t)$ converges
uniformly to $[\hbox{\rm Crit}(G^+,\Gamma_t)]$ as $j,k-j\to\infty$; that is
if $\psi$ is any test function and $\kappa_1,\kappa_2\to\infty$, then
$$\lim_{k\to\infty}\max_{\kappa_1(k)\le j\le k-\kappa_2(k)}\max_{t\in T}
\Big|\int\psi (\lambda^u_{j,k}(\beta,t)-[\hbox{\rm
Crit}(G^+,\Gamma_t)])\Big|=0.$$

\give Proof.  As $j,k-j\to\infty$, $\Gamma_t\cap
f^jZ_k(\alpha,\beta)$ converges to $\hbox{\rm Crit}(G^+,\Gamma_t)$ uniformly
in $\alpha\in\cV^u(j,t)$,
$t\in T$ and $\beta\in({\bf P}^1)^*$.  The
Lemma follows since $\cV^u(j,t)$ approaches full measure as $j\to\infty$. 
\qed
 
\proclaim Lemma 5.5.  Let $\Gamma$ be conformally equivalent to the unit disk
in $\cx{}$, and let $\Gamma'$ be a relatively compact open subset of
$\Gamma$.  Let $h:\Gamma\to\cx2$ be a holomorphic function such that
$$\max_{\zeta\in\Gamma}|h|\le C\min_{\zeta\in\Gamma}|h|$$
for some $C<\infty$.  Then there is a constant $0<b<1$, depending only on
$\Gamma'$, such that
$$m=\max_{\alpha\in{\bf
P}^1}\{\zeta\in\Gamma':h(\zeta)\cdot\alpha=0\}\eqno(5.1)$$
satisfies either $Cb^m\ge\sqrt 3/2$ or
$$\int_{\alpha\in{\bf P}^1} \sigma(\alpha)
\#\{\zeta\in\Gamma':h(\zeta)\cdot\alpha=0\}\le m\pi(2Cb^m)^2.\eqno(5.2)$$

\give Proof.  Without loss of generality, we may assume that
$\sup_\Gamma|h|=C$ and $\inf_\Gamma|h|=1$.  There is a number $0<b<1$
depending only on $\Gamma'$ such that for any holomorphic function $\psi$ on
$\Gamma$ with $m$ zeros in $\Gamma'$, 
$$\max_{\Gamma'}|\psi|\le b^m\max_{\Gamma}|\psi|.$$
Let us fix $\alpha_0$ with $|\alpha_0|=1$ such that the maximum is attained in
(5.1).  It follows that
$$\max_{\zeta\in\Gamma'}|h(\zeta)\cdot\alpha_0|\le b^mC.$$

If $\theta(\zeta)$ is the angle between $\hbox{\rm Ker}h(\zeta)$ and
$\alpha_0\in\cx2$, then 
$$|h|\sin\theta(\zeta)=|h(\zeta)\cdot\alpha_0|.$$
Since $|h(\zeta)|\ge1$, it follows that $|\sin\theta(\zeta)|\le Cb^m$.  It
follows that  $h(\zeta)\cdot\alpha\ne0$ for $\zeta\in\Gamma'$ if the sine of
the angle between $\alpha$ and $\alpha_0$ is greater than $Cb^m$.  If
$|\sin\theta|<\sqrt3/2$, then $\theta/2<|\sin\theta|$.  Thus if
$Cb^m<\sqrt3/2$, then $|\theta(\zeta)|\le2Cb^m$, and so 
$\alpha\mapsto\#\{\zeta\in\Gamma':h\cdot\alpha=0\}$ is supported in a disk of
radius $2Cb^m$ about $\alpha_0$.  In this case the integral in (5.2) is bounded
by
$m\pi(2Cb^m)^2$. \qed

\proclaim Lemma 5.6.  Let $\Gamma'\subset\Gamma$, $h$, $b$, and $C$ be as in
Lemma 5.5.  If $\cV\subset{\bf P}^1$ is contained in a disk of radius
$\delta$, then 
$$\int_{\alpha\in\cV} \sigma(\alpha)
\#\{\zeta\in\Gamma':h(\zeta)\cdot\alpha=0\}\le
\delta^2 C'\log({1\over\delta}),$$
where $C'$ depends only on $b$ and $C$.

\give Proof.  Let us choose $\alpha_0$ which maximizes $m$ in (5.1).  If
$Cb^m<\sqrt3/2$, then 
$$\int_{\alpha\in\cV} \sigma(\alpha)
\#\{\zeta\in\Gamma':h(\zeta)\cdot\alpha=0\}\le
m\pi\delta^2<{\pi\log(\sqrt3/2C)\over\log b}\delta^2.$$
If $Cb^m\ge\sqrt3/2$, then by Lemma 5.5 the integral is bounded by
$m\pi(2Cb^m)^2$.  We also have the trivial upper bound $m\pi\delta^2$.  Thus
$$\int_{\alpha\in\cV} \sigma(\alpha)
\#\{\zeta\in\Gamma':h(\zeta)\cdot\alpha=0\}\le\min(m\pi(2Cb^m)^2,m\pi\delta^2)
={\log(\delta/2C)\pi\delta^2\over\log
b}$$
since the minimum is attained when
$2Cb^m=\delta$.  \qed

\proclaim Lemma 5.7.  For $\kappa_0$ sufficiently large, there exists a
constant $C$ such that for $k-j\ge \kappa_0$, $h=(\beta\circ
f^{k-j})^{-1}\beta\circ Df^{k-j}:\cx2\to\cx2$ satisfies
$$\max_{\Gamma_t}|h|\le C\min_{\Gamma_t}|h|$$
for all $t\in T$ and $\beta\in({\bf P}^1)^*$.

\give Proof.  This is a direct consequence of Lemma 5.2. \qed

\proclaim Lemma 5.8.  
$$\lim_{{j\to\infty}\atop{k-j\to\infty}}\lambda^s_{j,k}(\beta,t)=0.$$

\give Proof.  Let $\Gamma'_t\subset\Gamma_t$ be a relatively compact open subset
with no critical points in the boundary.  By Lemma 5.5, 5.6, and 5.7, we have
$$\lambda^s_{j,k}(\beta,t)\contract\Gamma'_t\le  \max\left({1\over
b}\left(\log{\pi\over
4C}\right)\pi\delta_j^2,{\log(\delta_j(4C)^{-1})\over\log
b}\pi\delta_j^2\right)$$
where $\delta_j$ is chosen so that $\cV^s_j(\beta,t)$ is contained in a disk
of radius $\delta_j$.  The Lemma now follows since $\delta_j\to0$ as
$j\to\infty$ and since $\Gamma'_t$ can be chosen to exhaust $\Gamma_t$.  \qed

Define
$$\hat Z_k(\beta)=\int_{\alpha\in{\bf
P}^1}\sigma(\alpha)\,[Z_k(\alpha,\beta)].$$

\proclaim Theorem 5.9.  As $j,k-j\to\infty$ the restrictions of
$\mu^-\wedge f^j_*\hat Z_{k}(\beta)$ to $U^+$ 
converge to $\mu^-_c$ in the sense
of currents on $U^+$.

\give Proof.  We let $P^u=P$ be an unstable box as above.  We choose an unstable
box $P'\subset P$ such that $\Gamma'_t$ is relative compact in $\Gamma_t$,
and there are no critical points on $\partial\Gamma'_t$.  Further, we may
assume that $\mu^-_c(\bigcup_{t\in T}\partial\Gamma_t)=0$.  By Lemma 3.2 it
suffices to show that
$$\lim_{{j\to\infty}\atop{k-j\to\infty}}(\mu^-\wedge f^j_*\hat
Z_k(\beta))\contract P =\mu^-_c\contract P.$$
Using the notation above
$$\eqalign{
\mu^-\wedge f^j_*\hat Z_k(\beta)\contract P &=(\mu^-\contract P)\wedge
f^j_*\hat Z_k(\beta)\cr
&=\int_{t\in T}\mu^-_\tau(t)\,[\Gamma_t]\wedge f^j_*\hat Z_k(\beta)\cr
&=\int_{t\in T}\mu^-_\tau(t)\int_{\alpha\in{\bf
P}^1}\sigma(\alpha)\,[\Gamma_t]\wedge f^j_*[Z_k(\alpha,\beta)]\cr
&=\int_{t\in T}\mu^-_\tau(t)\int_{\alpha\in{\bf
P}^1}\sigma(\alpha)\,[\Gamma_t\cap f^jZ_k(\alpha,\beta)].\cr}$$
If we break up the inner integral as ${\bf
P}^1=\cV^s_j(\beta,t)\cup\cV^u_j(\beta,t)$, then we have
$$\mu^-\wedge f^j_*\hat Z_k(\beta)\contract P =\int_{t\in
T}\mu^-_\tau(t)\lambda^s_{j,k}(\beta,t) + \int_{t\in
T}\mu^-_\tau(t)\lambda^u_{j,k}(\beta,t).$$
It follows from Lemma 5.8, then, that the first integral on the right hand side
converges to zero, and from Lemma 5.4 that the second integral converges to
$\mu^-_c\contract P$.   \qed

We observe that as in (4.5) $(Df^k)^*\Theta=\int\sigma(\beta)\hat Z_k(\beta)$,
we may integrate the previous result with respect to $\beta$ to obtain:

\proclaim Corollary 5.10.  Let $\Theta$ be as in (4.4).  Then as
$j,k-j\to\infty$ the restrictions of the currents $\mu^-\wedge
f^j_*(Df^{k})^*\Theta$ to $U^+$ converge to $\mu^-_c$ in the sense of
currents on $U^+$.

\section 6. The Integral Formula

The main goal of this Section is to prove Theorem 6.1, which gives the formula
(0.2).  In fact, Theorem 6.1 is a consequence of Theorem 6.2, relating the rate
of expansion from (4.4) to the unstable critical measure.  This may be viewed as
applying Corollary 5.10 inside the integral formula of Proposition 4.6.

For a set $P$, we put $\tilde P=\bigcup_{n\in{\bf Z}}f^nP$.
We will say that a Borel set $P$ is a {\it fundamental domain }
for $\cC^u$ if
$\tilde P\supset\cC^u$ and if $P\cap f^nP=\emptyset$ for all $n\ne0$.

\proclaim Theorem 6.1.  Let $P\subset \cC^u$ be a
fundamental domain for $\cC^u$.  Then
$$\lambda^+(\mu)=\log d+\int_PG^+\mu^-_c.$$

\give Remark. A convenient choice for fundamental domain is $\{1\le
G^+<d\}\cap\cC^u$.  This
gives
$$\lambda^+(\mu)=\log d+\int_{\{1\le
G^+<d\}}G^+\mu^-_c.\eqno(6.1)$$

For a
domain $P$ satisfying $P\cap f^nP=\emptyset$ for all $n\ne0$, every
point $x\in\tilde P$ may be written uniquely as $x=f^ny$, so we have
a projection $\pi_P:\tilde P\to P$ given by $\pi_P(x)=y$; it is evident
that $\pi_P$ is Borel measurable.

\proclaim Theorem 6.2.  Let $P\subset J^-$ be a Borel set such that
$\mu^-_c(\partial P)=0$, where $\partial P$ denotes the boundary relative
to $J^-$.
If $P\cap f^nP=\emptyset$ for all $n\ne0$, then $$G^+\mu^-_c\contract P=
\lim_{k\to\infty}(\pi_P)_*\left({1\over k}G^+\mu^-\wedge(Df^k)^*\Theta
\contract\tilde P\right).\eqno(6.2)$$

\give Remark.  Both sides of the equation put no mass on $J^-\cap K$, so
without loss of generality we may assume that $P\subset J^-- K$.
Indeed,
the general case follows from the case where $P$ is a fundamental domain.

\give Proof of Theorem 6.1.  We will show how Theorem 6.1 is deduced from
Theorem 6.2.  We first prove
$$\lambda^+(\mu)=\log d + \int_{\{t\le G^+<td\}} G^+\mu^-_c$$
for some value of $t$.  For this, we note
that for every $t>0$,
$P_t=\{t\le  G^+<td\}\cap J^-$
is a fundamental domain for $J^--K$.
By the fact that $G^+$ is pluriharmonic, we have
$\partial P_t=\{G^+=t\}\cup\{G^+=td\}$,
so that the boundaries $\partial P_t$ are disjoint for $0<t<d$.
Now since $\mu^-_c$ is $\sigma$-finite,
we have $\mu^-_c(\partial P_t)=0$ for all but
countably many values of $t$.
So we may apply Proposition 4.6 and
Theorem 6.2 to conclude that the formula above holds for such $t$.

Now we conclude with the observation that if the Theorem holds for one
choice of Borel measurable fundamental domain, it holds for any other.
Given the fundamental domain $P$, the restriction of the mapping $\pi_P:\{t\le
G^+<td\}\cap\cC^u\to P$ is one to one and onto.  Since $G^+\mu^-_c$ is
$f$-invariant, it follows that it is invariant under $\pi_P$, and thus
$$\int_PG^+\mu^-_c=\int_{\{t\le G^+<t d\}}G^+\mu^-_c,$$
which completes the proof. \qed

By (1.2), $\beta\cdot Df^k(\alpha)=\beta_2\alpha_2d^ky^{d^k-1}+\cdots$, so
that if $\alpha_2\beta_2\ne0$, then the total mass of the intersection current
is
$$\int\mu^-\wedge[Z_k(\alpha,\beta)]=d^k-1.\eqno(6.3)$$
\proclaim Lemma 6.3.  If $\kappa_2(k)$ satisfies
$\lim_{k\to\infty}(\kappa_2(k)-\log_dk)=-\infty$,
then
$$\lim_{k\to\infty}{1\over k}\int_{\{G^+<d^{-k+\kappa_2}\}}
G^+\mu^-\wedge(Df^k)^*\Theta=0.$$

\give Proof.  By (6.3), the total mass of $\mu^-\wedge[Z_k(\alpha,\beta)]$ is
$d^k-1$ for almost every  $\alpha,\beta\in\cx2-\{0\}$.  Thus
$${1\over k}\int_{\{G^+<d^{-k+\kappa_2}\}} G^+\mu^-\wedge[Z_k(\alpha,\beta)]\le
{1\over k}d^{-k+\kappa_2} d^k$$
so the Lemma follows from the condition on $\kappa_2$ after integrating with
respect to $\alpha$ and $\beta$.  \qed

For a tangent vector $\alpha\in\cx2$ we define
$$Z_\infty(\alpha)=U^+\cap\{\partial G^+\cdot\alpha=0\}.$$
We note that since $Df^k$ and $\partial G^+$ are nonsingular
$Z_k(\alpha',\beta)\cap Z_k(\alpha'',\beta)=\emptyset$ and
$Z_\infty(\alpha')\cap Z_\infty(\alpha'')=\emptyset$ for all $\alpha',\alpha''$
which define distinct elements of ${\bf P}^1$.

\proclaim Lemma 6.4.  For each nonzero $\alpha\in\cx2$ the currents
$[Z_k(\alpha,\beta)]$ converge to $[Z_\infty(\alpha)]$ as currents on $U^+$,
uniformly in $\beta$.  That is, if $\psi$ is a test form with compact support
in $U^+$, then
$$\lim_{j\to\infty}\max_{\beta}\Big|\int\psi\wedge
([Z_k(\alpha,\beta)]-[Z_\infty(\alpha)])\Big|=0.$$

\give Proof. Since by Lemma 5.2 the projective images of the defining functions
of $Z_k(\alpha,\beta)$ converge uniformly, this gives the uniform convergence
of the currents.  
\qed

Let $V^+(R)=\{|y|>|x|,|y|>R\}$.  Since $G^+(x,y)=\log|y|+O(|y|^{-1})$ on
$V^+(R)$, it follows that
$$\partial G^+\cdot\alpha={\alpha_2\over y}+O(|y|^{-2}).\eqno(6.4)$$
Multiplying this by $y^2$, we see that for $R$ large and $\alpha=(1,\alpha_2)$
$$V^+(R)\cap Z_\infty(\alpha)=\{\alpha_2y+A_1(x,y)+\alpha_2 A_2(x,y)=0\}$$
where $A_1,A_2$ are bounded and holomorphic in $V^+(R)$.  Thus we have
$|dy/dx|\le c|y^{-1}|$ on $V^+(R)\cap Z_\infty(\alpha)$, and for $|\alpha_2|$
sufficiently small $V^+(R)\cap Z_\infty(\alpha)$ is a complex disk
$\{y=\varphi_\alpha(x):x\in D_\alpha\}$ satisfying
$${c'\over|\alpha_2|}\le
|\varphi_\alpha(x)|\le{c''\over|\alpha_2|}\eqno(6.5)$$

\proclaim Lemma 6.5.  For any $c>0$
$$\lim_{k\to\infty}\int_{\{G^+>c\}} G^+\mu^-\wedge
(Df^k)^*\Theta=\int_{\alpha\in{\bf P}^1}
\int_{\{G^+>c\}}G^+\mu^-\wedge[Z_\infty(\alpha)]<\infty.$$
 
\give Proof.  We will first show that for every $\beta$
$$ \lim_{k\to\infty}\int_{\alpha\in{\bf P}^1}\sigma(\alpha)
\int_{\{G^+>c\}}G^+\mu^-\wedge[Z_k(\alpha,\beta)] 
=\int_{\alpha\in{\bf
P}^1}\int_{\{G^+>c\}}G^+\mu^-\wedge[Z_\infty(\alpha)]<\infty.$$
Let us consider the regions
$\{G^+>c\}\cap\{|y|\le R\}$ and  $\{G^+>c\}\cap V^+(R)$ separately.  The
currents
$[Z_\infty(\alpha)]$ put no mass on $\{G^+=c\}\cup\{|y|=R\}$. 
Thus by Lemma 6.4, the integrals over the first region converge to the desired
limit as $k\to\infty$.

For the second region, we first check that the integral on the right hand
side is finite.  For $R$ large, $Z_\infty(\alpha)\cap V^+(R)$ is a complex
disk as in (6.5).  Thus $[Z_\infty(\alpha)\cap V^+(R)]$ has total mass 1. 
Thus for $\alpha=(1,\alpha_2)$ with $|\alpha_2|\le\epsilon$, the integral over
the second region is no larger than 
$$\int_{|\alpha_2|<\epsilon}
\log\left({c\over|\alpha_2|}\right)\sigma(\alpha)<\infty.$$

The convergence of the integrals holds because the disks
$Z_k(\alpha,\beta)\cap V^+(R)$ are close to the disks $Z_\infty(\alpha)\cap
V^+(R)$ throughout
$V^+(R)$, uniformly in $k$.  Since this convergence as $k\to\infty$ holds
uniformly in $\beta$, we may integrate with respect to $\beta$ to complete
the proof of the Lemma. \qed

\give Proof of Theorem 6.2.  Let us choose $P$ to be an unstable box for
which $\mu^-_c(\partial'P)=0$, and let us write
$\lambda_k=G^+\mu^-\wedge(Df^k)^*\Theta\contract\tilde P$.  Since we may
exhaust the $P$ in the hypothesis of the Theorem by a countable family of
such stable boxes, it suffices to show that
$$\lim_{k\to\infty}{1\over k}(\pi_P)_*\lambda_k=G^+\mu^-_c\contract P.$$
We may choose $\kappa_2(k)$ as in Lemma 6.3 so that
$$\lim_{k\to\infty}{1\over
k}(\pi_P)_*(\lambda_k\contract\{G^+<d^{-k+\kappa_2}\})=0.$$

Now for any positive integer $\kappa_1$ we set $c=d^{-\kappa_1}$, so by Lemma
6.5 the integrals $\int G^+\mu^-\wedge(Df^k)^*\Theta$, with $k\ge0$, are all
bounded by a number $m(\kappa_1)$.  We may define a function $\kappa_1(k)$ to
increase to infinity sufficiently slowly that $k^{-1}m(\kappa_1(k))\to0$ as
$k\to\infty$.  It follows, then, that
$$\lim_{k\to\infty}{1\over
k}(\pi_P)_*(\lambda_k\contract\{G^+>d^{-\kappa_1}\})=0.$$

Choosing $\kappa_1$ possibly smaller, we also have
$\lim_{k\to\infty}k^{-1}(\kappa_2+\kappa_1)=0$.  Now for $j=j_k$ satisfying
$\kappa_1\le j\le k-\kappa_2$ it follows from Theorem 5.9 that
$$\lim_{k\to\infty}f^j_*(\lambda_k\contract f^{-j}P)=G^+\mu^-_c\contract P.$$
Thus from the uniformity of the convergence in Theorem 5.9 we have
$$\eqalign{\lim_{k\to\infty}{1\over
k}(\pi_P)_*&(\lambda_k\contract\{d^{-k+\kappa_2}\le G^+\le d^{-\kappa_1}\}) \cr
&=\lim_{k\to\infty}{1\over k}\sum_{j=\kappa_1}^{k-\kappa_2}
f^j_*(\lambda_k\contract f^{-j}P)
=G^+\mu^-_c\contract P,\cr}$$
which completes the proof.  \qed

\proclaim Corollary 6.6. 
$$\lambda^-(f)=-\log d-\int_{\{1\le G^-<d\}}G^-\mu^+_c.$$

\give Proof. We can apply the integral formula to $f^{-1}$. The corresponding
invariant measure is the same, i.e.\ $\mu_f=\mu_{f^{-1}}$. Replacing $f$ by
$f^{-1}$ interchanges the role of stable and unstable directions and changes the
signs of the exponents. If we write
$\lambda^+(f)$ and
$\lambda^-(f)$ for the Lyapunov exponents of $f$ then we observe that
$\lambda^-(f)=-\lambda^+(f^{-1})$. Thus the integral formula applied to
$f^{-1}$ yields the above formula.

The following characterization is a consequence of the integral formula:

\proclaim Corollary 6.7.  The following are equivalent:
\item{1.}  $\lambda^+(\mu)=\log d$.
\item{2.}  $\mu^-_c=0$.
\item{3.}  For $\mu$ a.e.\ $x$, $G^+|_{W^u(x)-K^+}$ has no critical
points.

\give Proof.  The measure $\mu^-_c$ has all of its mass on the set $G^+>0$, so
by (6.1), if $\Lambda=\log d$, then $\mu^-_c=0$.  Thus (1) implies (2).  The
construction of the measure shows that if the measure vanishes, then
$G^+|_{W^u(x)-K^+}$ can have no critical points for $\mu$ a.e.\ $x$.  So (2)
implies (3).  Similarly, if $G^+|_{W^u(x)-K^+}$ has no critical points, then
the measure $\mu^-_c$ is zero, so (3) implies (1) by formula (6.1). \qed

Applying Corollary 6.7 to $f^{-1}$ gives: 

\proclaim Corollary 6.8.  The following are equivalent:
\item{1.}  $\lambda^-(\mu)=-\log d$.
\item{2.}  $\mu^+_c=0$.
\item{3.}  For $\mu$ a.e.\ $x$, $G^-|_{W^s(x)-K^-}$ has no critical
points.

In [BS6]  we will explore further the topological consequences of the
nonexistence of critical points. In particular we will show that if
$\mu^{\pm}_c=0$, then $\cW^u\cap U^+$ is in fact a locally trivial lamination, so
that the critical points satisfy $\cC^{s/u}=\emptyset$ in a strong
pointwise sense (not just on  unstable manifolds of Pesin regular points).

We close by noting some relations between the existence of stable and unstable
critical points and the Jacobian determinant of $f$. Recall that the Jacobian
determinant of a polynomial diffeomorphism, $\det Df_p$ is depends only on $f$
and not on the point $p$. If $|\det Df|<1$ we say that $f$ is dissipative. If
$|\det Df|=1$ we say that $f$ is volume preserving.

\proclaim Proposition 6.9. If $f$ is dissipative, then $\cC^s\ne\emptyset$. If
$f$ is volume preserving, then $\cC^s=\emptyset$ if and only if
$\cC^u=\emptyset$.

\give Remark. When $\det Df=1$ then $f$ is conjugate to its inverse so that the
equivalence of the conditions $\mu^+_c=0$ and $\mu^-_c=0$ is clear. 

\give Proof. It is a general property of Lyapunov exponents that the sum of the
exponents is related to the Jacobian determinant. We have
$\lambda^+(\mu)+\lambda^-(\mu)=\int\log|\det f|d\mu=\log|\det f|$. Combining
this fact with the integral formulas for
$\lambda^\pm(\mu)$ gives:
$$\int_{\{1\le G^+<d\}}G^+\mu^-_c -\int_{\{1\le G^-<d\}}G^-\mu^+_c =\log|\det
Df|.$$
The contribution of each integral is non-negative and is positive when the
corresponding measure is non-zero. When $f$ is dissipative the right hand side
of the equation is negative. It follows that the value of the second integral
must be non-zero hence the equivalent conditions of the Corollary 6.8 are all
false. When
$f$ is volume preserving then the right hand side of the equation is zero. It
follows that the value of the integrals are equal. Hence the equivalent
conditions of Corollary 6.7 are equivalent to those of Corollary 6.8.
\qed

\section  Appendix A: Lyapunov Exponent of Real Horseshoes

Let $f_R$ be a polynomial automorphism of degree $d$ with real coefficients,
so $f_R:{\bf R}^2\to{\bf R}^2$ has a real polynomial inverse.   Let us suppose
that there is a topological square $D\subset{\bf R}^2$ such that $f_R$
maps $D$ across itself $d$ times.  A heuristic version of the case $d=3$ is
shown in Figure 1; the horizontal lines represent stable manifolds.
This situation occurs for the mapping
$$f:(x,y)\mapsto(y,y^d+c_{d-2}y^{d-2}+\dots+c_0-ax)$$
in a non-empty real parameter region, for instance, if $d=2$ and $-c_0\gg0$
or if $d=3$ and $-c_1\gg|c_0|^{2/3}$ In this case $f$ has a weak $d$-fold
horseshoe, and if follows (see Friedland and Milnor [FM]) that $f_R$ is
topologically conjugate on the set
$K_R:=\bigcap_{n\in{\bf Z}}f^nB$ to the bilateral shift on $d$ symbols.  In
this case $f_R$ has topological entropy equal to $\log d$, and by [BLS1]
$K_R=J_C=K_C\subset{\bf R}^2$, where $J_C$ and $K_C$ denote the sets
$J$ and $K$ for the complexified mapping $f_C:\cx2\to\cx2$.

We let $V_1,\dots,V_{d}$ denote the (vertical) components of $D\cap fD$. 
Then there are components $B_1,\dots,B_{d-1}$ of $fD-J^+$ with the property that
$B_j$ intersects two distinct vertical components.  These are the fundamental
bends; the case $d=3$ is depicted in Figure 1.  We let $\cC_{0,j}$ denote the set
of critical points lying in the
$j$-th fundamental bend, i.e.\ $\cC_{0,j}=B_j\cap\cC^u$. 
Thus $\cC_0:=\cC_{0,1}\cup\cdots\cup\cC_{0,d-1}$ are all the critical points that
lie in the fundamental bends.  The critical points of the $n$-th image under $f$,
$n\in{\bf Z}$, are defined as $\cC_n=f^n\cC_0$.

\proclaim Lemma A.1.  Let $f$ be a $d$-fold real horsehoe, as above.  For
every $x\in J$, the restriction $G^+|W^u(x)$ has the property that every
component of $\{G^+|W^u(x)<c\}$ is relatively compact in $W^u(x)$. 

\give Proof.  Let $W^u(p)$ be the unstable manifold of a periodic point $p$. 
Since $\lim_{\zeta\to J}G^+(\zeta)=0$, and since $W^u(p)\cap J$ is a Cantor
set, there is a $\lambda_0>0$ such that the component of $\omega_0$ of
$W^u(p)\cap\{G^+<\lambda_0\}$ containing $p$ is relatively compact.  If
$\omega$ is any component of $W^u(p)\cap\{G^+<\lambda\}$, then
$f^{-n}\omega\subset\omega_0$ for $n$ sufficiently large.  Thus $\omega$ is
relatively compact.  Now the result follows since $f$ is hyperbolic, and $J$
has a local product structure.\qed

\vskip3.5in
\special {picture Figure1 scaled 800}

\proclaim Lemma A.2.  Let $E$ be a closed subset of $\cx{}$, and let
$L_j\subset{\bf R}$ be disjoint, open intervals such that ${\bf R}-E=\bigcup
L_j$.  Let $h\ge0$ with $h(x+iy)=h(x-iy)$ be continuous on $\cx{}$ and harmonic
on $\cx{}-E$, and let $E=\{h=0\}$.  If each connected component of $\{h<c\}$
is bounded, then for each $j$ there exists a unique critical point $c_j\in L_j$.
Further, the $\{c_j\}$ are all of the critical points of $h$.

\give  Proof.  Let $\omega$ be a component of $\{h<\lambda\}$.  Since
$\omega$ is relatively compact, it follows from the maximum principle that
$E\cap\omega\ne\emptyset$.  Since $\tilde\omega=\{\bar z:a\in\omega\}$ is
also a component of $\{h<\lambda\}$, and since $\emptyset\ne\tilde\omega\cap
E=\omega\cap E\subset{\bf R}$, it follows that $\tilde\omega=\omega$.

Now we claim that $\omega\cap{\bf R}$ is an interval.  It is nonempty, and if
it contains two components, then by the fact that $\omega$ is connected and
$\omega=\tilde\omega$, we have that $\cx{}-\omega$ contains a compact component. 
But this contradicts the maximum principle since $h$ is subharmonic on $\cx{}$.

Next suppose that there is a critical point $c\notin{\bf R}$.  Let
$\omega',\omega''$ be two components of $\{h<h(c)\}$ which contain $c$ in
their boundaries (possibly $\omega'=\omega''$).  Since these sets are
invariant under complex conjugation, it follows that $\bar c$ is also in their
boundaries.  Thus the complement of $\omega'\cup\omega''\cup\{c,\bar c\}$
in $\cx{}$ contains a compact component, which violates the maximum
principle.  Thus all critical points are real.

Let us fix an interval $L_j=(a_j,b_j)$ and let $\omega_{a_j}(\lambda)$
(resp.\ $\omega_{b_j}(\lambda)$) denote the component of $\{h<\lambda\}$
containing $a_j$ (resp.\ $b_j$).  For $\lambda>0$ sufficiently small,
$\bar\omega_{a_j}(\lambda)\cap \bar\omega_{b_j}(\lambda)=\emptyset$.  If
there is a critical point $c\in\bar\omega_{b_j}(\lambda)\cap L_j$, then we
may decrease $\lambda$ so that $c\in\partial\omega_{b_j}(\lambda)\cap L_j$. 
Since $\omega_{b_j}(\lambda)=\tilde\omega_{b_j}(\lambda)$, there must be a
distinct component $\omega$ of $\{h<\lambda\}$ such that $c\in\partial\omega$. 
But since $\omega\cap{\bf R}$ is an interval, and $\omega\cap E\ne\emptyset$, we
have  $a_j\in\omega$.  This is a contradiction, so we have no critical points in
$L_j\cap(\omega_{a_j}(\lambda)\cup\omega_{b_j}(\lambda))$ if
$\bar\omega_{a_j}(\lambda)\cap\bar\omega_{b_j}(\lambda)=\emptyset$.

On the other hand, if $\lambda$ is the supremum of numbers such that
$\bar\omega_{a_j}(\lambda)\cap\bar\omega_{b_j}(\lambda)=\emptyset$, then
$\bar\omega_{a_j}(\lambda)\cap\bar\omega_{b_j}(\lambda)=\{c_j\}$ is the unique
critical point in $L_j$.  \qed

The following theorem shows that the critial points for a horseshoe are
arranged in the fashion given schematically in Figure 2.

\vskip3.8in
\hskip.5in\special {picture Figure2d scaled 800}

\proclaim Theorem A.3.  Let $x\in J$ be given,
and let $W^u_R(x)$ denote the (real) unstable manifold passing through $x$. 
For each connected component $\gamma$ of $W^u_R(x)\cap B_j$ there is a unique
critical point $c_\gamma$ for the complexified mapping $f_C$.  The union of
all such critical points gives $\cC_0$, and $\cC^u=\bigcup_{n\in{\bf
Z}}\cC_n$.  In particular, all complex critical points are real.

\give Proof.  Let $W^u(x)$ denote the complex stable manifold throught $x$,
and let $\psi:\cx{}\to W^u(x)$ denote a uniformization.  Since $f$ is real,
we may replace $\psi(\zeta)$ by $\psi(e^{i\theta}\zeta)$ so that
$\psi_{\bf R}:{\bf R}\to W^u_R(x)$.  Let $h=G^+\circ\psi$
and let
$E=\psi^{-1}(W^u(x)\cap J)$ so that ${\bf R}-E=\bigcup L_j$.  Then $h$ is a
subharmonic function on $\cx{}$, and by Lemma A.1 it satisfies the hypotheses
of Lemma A.2.  It follows that all of the critical points of $h$ are real,
and so $\cC^u\cap W^u(x)\subset{\bf R}$.  Thus $\cC^u=\bigcup_{x\in
J}\cC^u\cap W^u(x)\subset{\bf R}$.  Also by Lemma A.2, we have that each
critical point $c\in\cW^u(x)$ corresponds uniquely to an interval $L_j$, and
$\psi(L_j)$ corresponds to a connected component $\gamma_c$ of
$W^u_R(x)-\cW^s=W^u_R(x)- J$.  Now it is a property of the horseshoe that
for each component $\gamma_c$, there is an $n\in{\bf Z}$ such that
$f^n\gamma_c\subset B_j$ for some $j$.  \qed

\give Remark.  If $f_R$ has the form above, then the line $\{x=0\}$ will
intersect the image of any non-horizontal line exactly once.  Under
iteration, this yields that $\{x=0\}$ will intersect each component of an
unstable manifold in a bend exactly once.  Since the total mass of the
intersection measure
$\mu^-\wedge[\{x=0\}]$ is 1, we see by Theorem A.1 that
$\mu^-_c(\cC_{0,j})=1$.  Further, $\mu^-_c$ has a balanced property that
allows us to define it in terms of the ``generational'' structure.  It
suffices to define $\mu^-_c$ on $\cC_{0,j}$, i.e.\ inside one of the fundamental
bends.  For this, we note that $B_j\cap f^nD$ has $d^{n-1}$ connected
components.  The intersection of any of these components with $\cC_0\cap B_j$
has mass $d^{-n+1}$, and this defines
$\mu^-_c$ on all Borel subsets of $\cC_0\cap B_j$.

\proclaim Theorem A.4.  If $f$ is a real horseshoe mapping as above, then the
Lyapunov exponent is given by
$$\Lambda = \log d + \int_{\cC_0} G^+\mu^-_c.$$
Further, we have the estimate
$$(d-1)\min_{\cC_0}G^+ < \Lambda - \log d <(d-1)\max_{\cC_0}G^+.$$

\give Proof.  The integral formula follows from Theorem 6.1 and Theorem A.1. 
The inequalities follow since $\mu^-_c(\cC_0))=d-1$, as was remarked above. 
\qed

\section Appendix B: Heteroclinic tangencies in $U^+\cap U^-$

We discuss the behavior of $f$ on $U^+\cap U^-$.  Conversations with J.H. Hubbard
have been helpful for our understanding of the critical locus in this region. 
The map $G^+:U^+\to{\bf R}^+$ has been studied as a fibration in [H, HO], where
it was shown that the level sets $\{G^+=c\}$ are foliated by complex manifolds
which are dense and conformally equivalent to $\cx{}$.  

It is shown in [H,HO] that we have an analytic function $\varphi^+$ on $V^+$
given by the formula
$$\varphi^+(x,y)=\lim_{n\to\infty}(\pi_2\circ f^n(x,y))^{1\over
d^n},$$ 
where we take the $d^n$-th root so that $\varphi^+(x,y)=y+o(1)$ holds on $V^+$.
It is immediate that $\varphi^+\circ f=(\varphi^+)^d$ and $\log|\varphi^+|=G^+$
hold on $V^+$.  In particular, $\varphi^+$ is locally constant on the leaves of
$\cG^+$. 

For $|\xi|>R$,
$$\Delta_\xi:=\{p\in V^+:\varphi^+(p)=\xi\}$$
is a complex disk, and  $f\Delta_\xi\subset \Delta_{\xi^d}$.  By the
trapping property of $V^+$ the global leaf
$L_\xi$ of $\cG^+$ which contains $\Delta_\xi$ has the form
$$L_\xi=\bigcup_{n=1}^\infty f^{-n}\Delta_{\xi^{d^n}},$$
and it is evident that $L_\xi\cap V^+=\bigcup L_{\xi'}$, where the union is
taken over all $\xi'$ such that $\xi'\xi^{-1}$ is a $d^n$th root of unity.

\proclaim Proposition B.1.  The global leaves of $\cG^+$ are the
super-stable manifolds of $f$.

\give Proof.  By Lemma 1.2, $Df^n|_{T\cG^+}$ decreases super-exponentially as
$n\to+\infty$.  Thus any two points $\zeta',\zeta''$ in the same global leaf
of $\cG^+$ approach each other super-exponentially as $n\to+\infty$. 
Conversely, suppose that $\zeta',\zeta''\in U^+$ are not in the same global
leaf $U^+$.  Then for $n\ge n_0$, $f^n\zeta',f^n\zeta''\in V^+$, but 
$$\varphi^+(f^n\zeta')= \varphi^+ (f^{n_0}\zeta')^{n-n_0}\ne \varphi^+
(f^{n_0}\zeta'')^{n-n_0} = \varphi^+(f^n\zeta'').$$
Thus $\varphi^+(f^n\zeta')$ does not tend to $\varphi^+(f^n\zeta'')$, and
since $\varphi^+\approx y$ on $V^+$, it follows that $f^n\zeta'$ does not
tend to $f^n\zeta''$. \qed

The 2-form $\partial G^+\wedge\partial G^-$ is invariant under $f$, and its zero
locus defines the dynamical critical locus of $f$:
$$\cC:=\{(x,y)\in U^+\cap U^-:\tau^+=\tau^-\}=\{\partial G^+\wedge\partial
G^-=0\}.$$ 
Thus the critical locus consists of the points where the forward and
backward critical directions coincide; thus it can be thought of as the set of
heteroclinic tangencies of the super-stable and super-unstable manifolds.  

For $\epsilon>0$ there exists $R_\epsilon$ such that
$$(G^+,G^-)\approx(\log|y|,\log|x|), {\rm\ \ and\ \ }\partial G^+\wedge\partial
G^-\approx{dx\wedge dy\over 4xy}\eqno(B.1)$$ 
for $\epsilon|x|<|y|<\epsilon^{-1}|x|,\
|x|>R_\epsilon$, and so 
$$\cC\cap\{\epsilon|x|<|y|<\epsilon^{-1}|x|,\
|x|>R_\epsilon\}=\emptyset.\eqno(B.2)$$

The inclusions of sets $\iota_\pm:U^+\cap U^-\to U^\pm$ induce mappings on
homology
$${\iota_\pm}_*:H_1(U^+\cap U^-,{\bf Z})\to H_1(U^\pm,{\bf Z}).\eqno(B.3)$$

\proclaim Lemma B.2.  The mapping (B.3) is surjective, and $H_1(U^+\cap
U^-,{\bf Z})$ is not finitely generated.

\give Proof.  For $R$ large, consider the curve
$\gamma:\theta\mapsto(Re^{i\theta},Re^{i\theta})$, which is contained in
$V^+\cap V^-\subset U^+\cap U^-$.  Now $\varphi^+_*\gamma$ is approximately
the circle of radius $R$ in $\cx{}-\bar\Delta$, so it defines a
nontrivial homology class, and thus $\gamma$ defines a nontrivial element of
$ H_1(U^\pm,{\bf Z})$ in the range of ${\iota_\pm}_*$.  Since the range is
nonzero and invariant under $f_*^k$ for
$k\in{\bf Z}$, the maps ${\iota_\pm}_*$ are onto.   Finally, $H_1(U^+\cap
U^-,{\bf Z}) $ is not finitely generated because its image is not. (See [HO]
for these last two facts about $H_1(U^\pm,{\bf Z})$.) \qed

\proclaim Proposition B.3.  $\cC\ne\emptyset$.

\give Proof.  We consider the fibration
$$G=(G^+,G^-):U^+\cap U^-\to{\bf R}^+\times {\bf R}^+,$$
which has compact fiber.  By (B.1), the fiber  of $G$ over points of
$\{\epsilon|x|<|y|<\epsilon^{-1}|x|,\ |x|>R_\epsilon\}$ is a 2-torus.  If
$\cC=\emptyset$, then $dG^+\wedge dG^-\ne0$, and  the fibration is
locally trivial.  Since the base of the fibration is topologically trivial, it
follows that
$$H_1(U^+\cap U^-,{\bf Z})\cong H_1({\bf T}^2,{\bf Z})\cong {\bf Z}^2.$$
But this is not possible since, by Lemma B.2, $H_1(U^+\cap U^-,{\bf Z})$ is not
finitely generated.  Thus $\cC\ne\emptyset$. \qed

\proclaim Proposition B.4.  $\bar\cC\cap J^+\cap U^-\ne\emptyset$ and
$\bar\cC\cap J^-\cap U^+\ne\emptyset$.

\give Proof. If $\bar\cC\cap
J^-\cap U^+=\emptyset$, then $\cC$ is a closed subvariety of $U^+$.  Let
$\cC'=\cC\cap\{|y|>R,\,|y|>|x|\}$.  By (B.2), $\pi_2|_{\cC'}:\cC'\to\{|y|>R\}$
is proper, so it has some degree $\delta$.   This degree multiplies by $d$
under the mapping $f$. But since $f\cC=\cC$ this degree must stay constant. 
Thus we conclude that $\bar\cC\cap
J^-\cap U^+\ne\emptyset$.  The argument to show $\bar\cC\cap J^+
\cap U^-\ne\emptyset$
is the same. \qed

\centerline{\bf References}

\item{[BS1]} E.\ Bedford and J.\ Smillie, Polynomial diffeomorphisms of $\cx2$:
Currents, equilibrium measure and hyperbolicity. Invent.\ Math. 87, 69--99
(1990).

\item{[BS3]} E.\ Bedford and J.\ Smillie, Polynomial diffeomorphims of $\cx 2$
III: Ergodicity, exponents and entropy of the equilibrium measure.  Math. Ann
294, 395--420 (1992).

\item{[BS6]} E. Bedford and J. Smillie, Polynomial diffeomorphisms of $\cx2$
VI: Connectivity of J. 

\item{[BS7]} E. Bedford and J. Smillie, Polynomial diffeomorphisms of $\cx2$
VII: Hyperbolicity and external rays.

\item{[BLS1]} E.\ Bedford, M.\ Lyubich, and J.\ Smillie, Polynomial
diffeomorphims of $\cx 2$. IV: The measure of maximal entropy and laminar
currents.  Invent. math. 112, 77--125 (1993).

\item{[BLS2]} E.\ Bedford, M.\ Lyubich, and J.\ Smillie,  Distribution of
periodic points of polynomial diffeomorphisms of $\cx 2$, Invent. math. 114,
277--288 (1993).

\item{[BC]}  M. Benedicks and L. Carleson,  The dynamics of the H\'enon map,
Ann. Math. 133, 73--169 (1991).

\item{[BY]} M. Benedicks and L-S Young,  Sinai-Bowen-Ruelle measures for
certain H\'enon maps.  Invent. math.  112, 541--576 (1993).

\item{[Br]}  H.\ Brolin,  Invariant sets under iteration of rational functions.
 Ark. Mat. 6, 103--144 (1965).

\item{[C]} E. Chirka,  {\sl Complex Analytic Sets}, Kluwer, 1985. 

\item{[DH]}  A. Douady and J.H. Hubbard, It\'eration des polyn\^omes
quadratiques complexes.  C.R. Acad. Sci. Paris S\'erie I 294, 123--126 (1982).

\item{[FM]} S.\ Friedland and J.\ Milnor, Dynamical properties of plane
polynomial automorphisms.
Ergodic Theory Dyn. Syst. 9, 67--99 (1989).

\item{[H]} J.H.\ Hubbard,  H\'enon mappings in the complex domain, in: Chaotic
Dynamics and Fractals, ed. M. Barnsley \& S. Demko, Academic Press,
pp.\ 101--111 (1986).

\item{[HO]} J.H.\ Hubbard and R.\ Oberste-Vorth, H\'enon mappings in the complex
domain I: The global topology of dynamical space, Inst.\ Hautes \'Etudes Sci.\
Publ.\ Math.\ 79, 5--46 (1994).

\item{[LY]}  F.\ Ledrappier \& L.--S.\ Young,
The metric entropy of diffeomorphisms, I \& II,
 Annals of Math. 122
 509--539 \& 540-574 (1985).

\item{[Mn]} A.\ Manning,  The dimension of the maximal measure for a polynomial
map, Ann. Math. 119, 425--430 (1984).

\item{[Po]}  M.\ Pollicott,  {\sl Lectures on ergodic theory and Pesin theory
on compact manifolds}, London Mathematical Society Lecture Note Series 180,
Cambridge U. Press, 1993.

\item{[Pr]} F.\ Przytycki, Hausdorff dimension of harmonic measure on the
boundary of an attractive basin for a holomorphic map. Invent. Math. 80,
161--179 (1985).

\item{[Si]}  N. Sibony, Iteration of polynomials, U.C.L.A. course lecture notes.

\item{[T]} P. Tortrat, Aspects potentialistes de l'it\'eration des
polyn\^omes.  In: S\'eminaire de Th\'eorie du Potentiel Paris, No.\ 8 (Lect.
Notes Math., vol.\ 1235) Berlin Heidelberg New York: Springer 1987.

\item{[W]} H. Wu,  Complex stable manifolds of holomorphic diffeomorphisms,
Indiana U. Math. J. 42, 1349--1358 (1993).

\bye